\documentclass[final]{siamltex}         

\usepackage[dvips]{color}

\usepackage{graphicx}

\usepackage{subfigure}

\graphicspath{{C:/NS/Pictures/}}

\title{The Cauchy problem for the 3D Navier - Stokes equations. \\ New approach to the solution and its justification.}  

\author{ A. Tsionskiy, M. Tsionskiy  }               

\newcommand{\sz}[1]{\mbox{\Large #1 \normalsize}}
\newcommand{\szp}[1]{\mbox{\large #1 \normalsize}}

\begin{document}             

\maketitle                   

\begin{abstract}
Some known results regarding the Euler and Navier-Stokes equations were obtained by different authors. 
Existence and smoothness of solutions for the Navier-Stokes equations in two dimensions have been known for a long time. 
Leray $\cite{jL34}$ showed that the Navier-Stokes equations in three space dimensions have a weak solution. 
Scheffer $\cite{vS76}, \cite{vS93}$ and Shnirelman $\cite{aS97}$ obtained weak solution of the Euler equations with compact support in spacetime. Caffarelli, Kohn and Nirenberg $\cite{CKN82}$ improved Scheffer's results, and F.-H. Lin $\cite{fL98}$ simplified the proof of the results of J. Leray. Many problems and conjectures about behavior of weak solutions of the Euler and Navier-Stokes equations are described in the books of Bertozzi and Majda $\cite{BM02}$, Constantin $\cite{pC01}$ or Lemarié-Rieusset $\cite{pL02}$.

Solutions of the Navier-Stokes and Euler equations with initial conditions (Cauchy problem) for 2D and 3D cases were obtained in the convergence series form by analytical iterative method using Fourier and Laplace transforms in paper $\cite{TT10}$. These solutions were received in a form of infinitely differentiable functions, and that allows us to analyze all aspects of the problem on a much deeper level and with more details. Also such smooth solutions satisfy the conditions required in $\cite{CF06}$ for the problem of Navier-Stokes equations. 

For several combinations of problem parameters numerical results were obtained and presented as graphs $\cite{TT10}$,$\;\cite{TT11}$.

This paper describes detailed proof of convergence of the analitical iterative method for solution of the Cauchy problem for the 3D Navier - Stokes equations. The convergence is shown for wide ranges of the problem's parameters. Estimated formula for the border of convergence area  of the iterative process in the space of system parameters is obtained. Also we have provided justification of the analytical iterative method solution for Cauchy problem for the 3D Navier-Stokes equations.

\end{abstract}



\pagestyle{myheadings}
\thispagestyle{plain}
\markboth{A. TSIONSKIY, M. TSIONSKIY}{NEW APPROACH. SOLUTION AND ITS JUSTIFICATION.}

\section{Introduction}\ 

Approach to the solution of Cauchy problem for the 3D Navier-Stokes equations described in this paper is based on the form of differential equations in the statement of the problem and also limitations for the initial conditions and applied force. It grows from classic definition of function and classic methods of analysis.

First we have moved non-linear parts of equations to the right sides and then solved Cauchy problem for the Navier-Stokes equations by analytical iterative method. By doing that we have obtained more and more precise results for the right (non-linear) parts of equations on each step of iterative process.

For this purpose we have solved the system of linear partial differential equations with constant coefficients on every step of iterative process. We have obtained the solution using Fourier transforms for the space coordinates and Laplace transform for time.
From theorems about application of Fourier and Laplace transforms for system of linear partial differential equations with constant coefficients we see that in case if initial conditions and applied force are smooth enough functions decreasing in infinity, then the solution of the system of linear partial differential equations with constant coefficients is also a smooth function. (Corresponding theorems are presented in S. Bochner $\cite{SB59}$, V.P. Palamodov $\cite{VP70}$, G.E. Shilov $\cite{gS01}$, L. Hormander $\cite{LH83}$, S. Mizohata $\cite{SM73}$, J.F. Treves $\cite{JFT61}$).

In the problem statement for Navier-Stokes equations initial conditions and applied force are infinitely differentiable functions decreasing rapidly to zero in infinity. Hence, all these theorems are applicable in this case.
Next, by using theorems from $\cite{gS01}$ regarding Fourier transforms for infinitely differentiable and decreasing to zero in infinity functions we have shown that all functions of Cauchy problem for Navier-Stokes equations are staying infinitely differentiable and decreasing to zero in infinity also rapidly enough in Fourier transforms on any step of iterative process.

After that we have received a superior estimation for the solution of the Cauchy problem for the 3D Navier - Stokes equations
by iterative method. The purposes of this estimation are:

$\;\;\;\;\;$1)	to show convergence of the iterative method;

$\;\;\;\;\;$2)	to obtain analytical form of the first and second steps of the iterative process;

$\;\;\;\;\;$3)	to receive estimated formula for the border of convergence region  of the iterative process in the

space of system parameters.

Also we have provided justification for the solution by analytical iterative method of Cauchy problem for the 3D Navier-Stokes equations.
While doing so, we have introduced perfect spaces of functions and vector-functions (I.M. Gel'fand, G.E. Chilov $\cite{GC68}$), in which we have looked for the solution of the problem.
We have demonstrated equivalence of the solution of Cauchy problem in forms of differential and integral equations.
As a further step we have proved existence and uniqueness of the solution of Cauchy problem in all time range [0,$\infty$) by using the fixed point principle ( L.V. Kantorovich, G.P. Akilov $\cite{KA64},\;$ V.A. Trenogin  $\cite{VT80},\;$W. Rudin $\cite{WR73},\;$W.A. Kirk and B. Sims$\; \cite{KS01},\;$A. Granas and J. Dugundji$\; \cite{GD03},\;$J.M. Ayerbe Toledano, T. Dominguez Benavides, G. Lopez Acedo$\;\cite{ADL97}\;$) .
For this purpose three following theorems were proven in this paper:

Theorem 1: Integral operator of the problem is a contraction operator;

Theorem 2: Existence and uniqueness of the solution of the problem is valid for any t $\in$ [0,$\infty$);

Theorem 3: Solution of the problem is depending on t continuously.

 By using the concept of comparable norms we have shown that the energy of the whole process has a finite value for any t $\in$ [0,$\infty$).

\section{The mathematical setup}\ 

The Navier-Stokes equations describe the motion of a fluid in \(R^{N}\; ( N = 3 ) \). We look for a viscous incompressible fluid filling all of \(R^{N}\)  here. The Navier-Stokes equations are then given by

\begin{equation}\label{eqn1}
\frac{\partial u_{k}}{\partial t} \; + \; \sum_{n=1}^{N} u_{n}\frac{\partial u_{k}}{\partial x_{n}}\; =\;\nu\Delta u_{k}\; - \; \frac{\partial p}{\partial x_{k}} \; + \;f_{k}(x,t)\;\;\;\;\;(x\in R^{N},\;\;t\geq 0,\;\;{1\leq k \leq N}) 
\end{equation}
\begin{equation}\label{eqn2}
\emph{div}\,\vec{u}\;=\; \sum_{n=1}^{N} \frac{\partial u_{n}}{\partial x_{n}}\; =\;0\;\;\;\;\;\;\;\;\;\;  (x\in R^{N},t\geq 0)
\end{equation}

with initial conditions

\begin{equation}\label{eqn3}
\vec{u}(x,0)\; = \; \vec{u}^{0}(x)\;\;\;\;\;\;\;\;\;\; (x\in R^{N})
\end{equation}

Here \(\vec{u}(x,t)=(u_{k}(x,t)) \in R^{N},\;\; ({1\leq k \leq N})  \;-\; \)is an unknown velocity vector  \(( N = 3 ),\; p\,(x,t)\;-\;\) is an unknown pressure, \(\vec{u}^{0}(x)\;\) is a given, \(C^{\infty}\) divergence-free vector field \(,\; f_{k}(x,t)\;\)are components of a given, externally applied force \(\vec{f}(x,t)\),  \(\nu\) is a positive coefficient of the viscosity (if \(\nu = 0\) then $(\ref{eqn1})$ - $(\ref{eqn3})$ are the Euler equations), and \( \Delta\;=\; \sum_{n=1}^{N} \frac{\partial^{2}}{\partial x_{n}^{2}}\;\) is the Laplacian in the space variables. Equation $(\ref{eqn1})$ is Newton's law for a fluid element subject. Equation $(\ref{eqn2})$ says that the fluid is incompressible. For physically reasonable solutions, we accept 

\begin{equation}\label{eqn4}
u_{k}(x,t) \rightarrow 0\;\;, \;\;
\frac{\partial u_{k}}{\partial x_{n}} \;\rightarrow\;0\;\; 
\rm {as} 
\;\;\mid x \mid \;\rightarrow\; \infty\;\;\;( {1\leq k \leq N} ,\;\; {1\leq n \leq N}) \;\;\; 
\end{equation}

Hence, we will restrict attention to initial conditions $\vec{u}^{0}$ and force $\vec{f}$ that satisfy

\begin{equation}\label{eqn5}
\mid\partial_{x}^{\alpha}\vec{u}^{0}(x)\mid\;\leq\;C_{\alpha K}(1+\mid x \mid)^{-K} \quad \rm{on }\;R^{N}\;\rm{ for\; any }\;\alpha \;\rm{ and \;}K.
\end{equation}
and
\begin{equation}\label{eqn6}
\mid\partial_{x}^{\alpha}\partial_{t}^{\beta}\vec{f}(x,t)\mid\;\leq\;C_{\alpha \beta K}(1+\mid x \mid +t)^{-K} \quad \rm{on }\;R^{N}\times[0,\infty)\; \rm{ for \;any }\;\alpha,\beta \;\rm{ and \;}K.
\end{equation}

We add (\( - \sum_{n=1}^{N} u_{n}\frac{\partial u_{k}}{\partial x_{n}}\;\)) to both sides of the equations (\ref{eqn1}). Then we have:

\begin{equation}\label{eqn7}
\frac{\partial u_{k}}{\partial t}\;=\;\nu\,\Delta\,u_{k}\;-\;\frac{\partial p}{\partial x_{k}}\;+\;f_{k}(x,t)- \; \sum_{n=1}^{N} u_{n}\frac{\partial u_{k}}{\partial x_{n}}\;\;\;\;\;\;\;\;\;\;\;\; (x\in R^{N},\;\;t\geq 0,\;\;{1\leq k \leq N})
\end{equation}
\begin{equation}\label{eqn8}
\emph{div}\,\vec{u}\;=\; \sum_{n=1}^{N} \frac{\partial u_{n}}{\partial x_{n}}\; =\;0\;\;\;\;\;\;\;\;\;\;  (x\in R^{N},t\geq 0)
\end{equation}
\begin{equation}\label{eqn9}
\vec{u}(x,0)\; = \; \vec{u}^{0}(x)\;\;\;\;\;\;\;\;\;\; (x\in R^{N})
\end{equation}

\begin{equation}\label{eqn10}
u_{k}(x,t) \rightarrow 0\;\;, \;\;\frac{\partial u_{k}}{\partial x_{n}}\;\rightarrow\;0\;\; \rm{as} \;\;\mid x \mid \;\rightarrow\; \infty\;\;\;( {1\leq k \leq N} ,\;\; {1\leq n \leq N}) \;\;\; 
\end{equation}

\begin{equation}\label{eqn11}
\mid\partial_{x}^{\alpha}\vec{u}^{0}(x)\mid\;\leq\;C_{\alpha K}(1+\mid x \mid)^{-K} \quad\rm{on }\;R^{N}\;\rm{ for \; any }\;\alpha\;\rm{ and }\;K.
\end{equation}

\begin{equation}\label{eqn12}
\mid\partial_{x}^{\alpha}\partial_{t}^{\beta}\vec{f}(x,t)\mid\;\leq\;C_{\alpha \beta K}(1+\mid x \mid +t)^{-K} \quad\rm{on }\;R^{N}\times[0,\infty)\;\rm{ for \; any }\;\alpha,\beta\;\rm{ and }\;K.
\end{equation}

We shall solve the system of equations (\ref{eqn7}) - (\ref{eqn12}) by the analitical iterative method. To do so we write this system of equations in the following form:

\begin{equation}\label{eqn13}
\frac{\partial u_{jk}}{\partial t}\;=\;\nu\,\Delta\,u_{jk}\;-\;\frac{\partial p_{j}}{\partial x_{k}}\;+\;f_{jk}(x,t) \;\;\;\;\;\;\;\;\;\;\;\; (x\in R^{N},\;\;t\geq 0,\;\;{1\leq k \leq N})
\end{equation}
\begin{equation}\label{eqn14}
\emph{div}\,\vec{u}_{j}\;=\; \sum_{n=1}^{N} \frac{\partial u_{jn}}{\partial x_{n}}\; =\;0\;\;\;\;\;\;\;\;\;\;  (x\in R^{N},t\geq 0)
\end{equation}
\begin{equation}\label{eqn15}
\vec{u}_{j}(x,0)\; = \; \vec{u}^{0}(x)\;\;\;\;\;\;\;\;\;\; (x\in R^{N})
\end{equation}

\begin{equation}\label{eqn16}
u_{j k}(x,t) \rightarrow 0\;\;, \;\;\frac{\partial u_{j k}}{\partial x_{n}}\;\rightarrow\;0\;\; \rm{as} \;\;\mid x \mid \;\rightarrow\; \infty\;\;\;( {1\leq k \leq N} ,\;\; {1\leq n \leq N}) \;\;\; 
\end{equation}

\begin{equation}\label{eqn17}
\mid\partial_{x}^{\alpha}\vec{u}^{0}(x)\mid\;\leq\;C_{\alpha K}(1+\mid x \mid)^{-K} \quad\rm{on }\;R^{N}\;\rm{ for \; any }\;\alpha\;and\;K.
\end{equation}

\begin{equation}\label{eqn18}
\mid\partial_{x}^{\alpha}\partial_{t}^{\beta}\vec{f}(x,t)\mid\;\leq\;C_{\alpha \beta K}(1+\mid x \mid +t)^{-K} \quad\rm{on }\;R^{N}\times[0,\infty)\;\rm{ for \; any }\;\alpha,\beta\;and\;K.
\end{equation}

Here j is the number of the iterative process step (j = 1,2,3,...).

\begin{equation}\label{eqn19}
f_{jk}(x,t)\; = \; f_{k}(x,t) \; - \; \sum_{n=1}^{N} u_{j-1,n}\frac{\partial u_{j-1,k}}{\partial x_{n}}\;\;\;\;\;\;({1\leq k \leq N})
\end{equation}

or the vector form

\begin{equation}\label{eqn20}
\vec{f}_{j}(x,t)\; = \; \vec{f}(x,t) \; - \;(\;\vec{u}_{j-1}\;\cdot\;\nabla\;)\;\vec{u}_{j-1}\;
\end{equation}
\\
For the first step of the iterative process (j = 1) we have:  

\[\\(\vec{u}_{0}\;\cdot\;\nabla\;)\;\vec{u}_{0}\;=\;0\]
and
\[\\\vec{f}_{1}(x,t)\; = \; \vec{f}(x,t) \]

\section{Solution and Estimation}\ 

Let us assume that all operations below are valid. We will prove the validity of these operations in next paragraphs.

We will solve the system of equations $(\ref{eqn13})\; - \;(\ref{eqn20})$ by analytical iterative method with the condition that $f_{jk}(x,t)$ is known
function on any step j of iterative process. Hence, we have that $(\ref{eqn13})\; - \;(\ref{eqn15})$ is the system of linear partial differential equations with constant coefficients. 

Because of that on the first stage we use Fourier transform $(\ref{A3})$ for solution of equations $(\ref{eqn13})\; - \;(\ref{eqn20})$ and get:

\[\\U_{jk}( \gamma_{1}, \gamma_{2}, \gamma_{3}, t)\;=\;F[u_{jk} ( x_{1}, x_{2}, x_{3}, t)]\]

\[F\sz{[} \frac{\partial^{2}u_{jk}( x_{1}, x_{2}, x_{3}, t)}{\partial x^{2}_{s}}\sz{]} \;=\;-\gamma^{2}_{s}U_{jk}( \gamma_{1}, \gamma_{2}, \gamma_{3}, t) \;\;\;\;\; \rm{[use (\ref{eqn16})]}\]

\[\\U_{k}^{0}( \gamma_{1} ,\gamma_{2} ,\gamma_{3})\;=\;F[u_{k}^{0} ( x_{1}, x_{2}, x_{3})]\]

\[\\P_{j}( \gamma_{1}, \gamma_{2}, \gamma_{3}, t)\;=\;F[p_{j}\, ( x_{1}, x_{2}, x_{3}, t)]\]

\[\\F_{jk}( \gamma_{1}, \gamma_{2}, \gamma_{3}, t)\;=\;F[f_{jk} ( x_{1}, x_{2}, x_{3}, t)]\]

\[ \\  k,s\;=\;1,2,3 \]

and then:

\begin{equation}\label{eqn134}
\frac{d U_{j1}( \gamma_{1}, \gamma_{2}, \gamma_{3}, t )}{d t}  \;=\;-\nu
( \gamma_{1}^{2} +\gamma_{2}^{2} +\gamma_{3}^{2}) U_{j1}( \gamma_{1}, \gamma_{2}, \gamma_{3}, t )\;+\;i\gamma_{1} P_{j}( \gamma_{1}, \gamma_{2}, \gamma_{3}, t )\;+\; F_{j1}( \gamma_{1}, \gamma_{2}, \gamma_{3},t )
\end{equation}
\begin{equation}\label{eqn135}
\frac{d U_{j2}( \gamma_{1}, \gamma_{2}, \gamma_{3}, t )}{d t}  \;=\;-\nu
( \gamma_{1}^{2} +\gamma_{2}^{2} +\gamma_{3}^{2}) U_{j2}( \gamma_{1}, \gamma_{2}, \gamma_{3}, t )\;+\;i\gamma_{2} P_{j}( \gamma_{1}, \gamma_{2}, \gamma_{3}, t )\;+\; F_{j2}( \gamma_{1}, \gamma_{2}, \gamma_{3},t )
\end{equation}
\begin{equation}\label{eqn136}
\frac{d U_{j3}( \gamma_{1}, \gamma_{2}, \gamma_{3}, t )}{d t}  \;=\;-\nu
( \gamma_{1}^{2} +\gamma_{2}^{2} +\gamma_{3}^{2}) U_{j3}( \gamma_{1}, \gamma_{2}, \gamma_{3}, t )\;+\;i\gamma_{3} P_{j}( \gamma_{1}, \gamma_{2}, \gamma_{3}, t )\;+\; F_{j3}( \gamma_{1}, \gamma_{2}, \gamma_{3},t )
\end{equation}

\begin{equation}\label{eqn137}
\gamma_{1} U_{j1}( \gamma_{1}, \gamma_{2}, \gamma_{3}, t ) \;+\; \gamma_{2}\, U_{j2}( \gamma_{1}, \gamma_{2}, \gamma_{3}, t ) \;+\; \gamma_{3}\, U_{j3}( \gamma_{1}, \gamma_{2}, \gamma_{3}, t ) \;=\;0
\end{equation}

\begin{equation}\label{eqn138}
U_{j1}(\gamma_{1}, \gamma_{2},  \gamma_{3},  0)\;=\; U_{1}^{0}(\gamma_{1} ,\gamma_{2} ,\gamma_{3})
\end{equation}

\begin{equation}\label{eqn139}
U_{j2}(\gamma_{1}, \gamma_{2},  \gamma_{3},  0)\;=\; U_{2}^{0}(\gamma_{1} ,\gamma_{2} ,\gamma_{3})
\end{equation}

\begin{equation}\label{eqn140}
U_{j3}(\gamma_{1}, \gamma_{2},  \gamma_{3},  0)\;=\; U_{3}^{0}(\gamma_{1} ,\gamma_{2} ,\gamma_{3})
\end{equation}

Hence we have received a system of linear ordinary differential equations with constant coefficients in regard to Fourier transforms $(\ref{eqn134})\;-\; (\ref{eqn140})\;$. At the same time the initial conditions are set only for Fourier transforms of velocity components $U_{j1}( \gamma_{1}, \gamma_{2}, \gamma_{3}, t ), \;U_{j2}( \gamma_{1}, \gamma_{2}, \gamma_{3}, t ), \;U_{j3}( \gamma_{1}, \gamma_{2}, \gamma_{3}, t )$. Because of that we can eliminate Fourier tranform for pressure $P_{j}( \gamma_{1}, \gamma_{2}, \gamma_{3}, t )$ from equations $(\ref{eqn134})\;-\; (\ref{eqn136})\;$ on the second stage of solution.
From here assuming that $\gamma_{1} \neq 0,\;\gamma_{2} \neq 0,\;\gamma_{3} \neq 0$, we eliminate $P_{j}( \gamma_{1}, \gamma_{2}, \gamma_{3}, t )$ from equations $(\ref{eqn134})\;-\; (\ref{eqn136})\;$ and find:

\begin{eqnarray}\label{eqn141}
\frac{d}{dt} \szp{[} U_{j2}( \gamma_{1}, \gamma_{2}, \gamma_{3}, t )\; -\;\frac{\gamma_{2}}{\gamma_{1}} \,U_{j1}( \gamma_{1}, \gamma_{2}, \gamma_{3}, t) \szp{]} \;= -\nu( \gamma_{1}^{2} +\gamma_{2}^{2} +\gamma_{3}^{2}) \szp{[} U_{j2}( \gamma_{1}, \gamma_{2}, \gamma_{3}, t )\; -\; \nonumber \\ 
\nonumber\\
-\;\frac{\gamma_{2}}{\gamma_{1}}\, U_{j1}( \gamma_{1}, \gamma_{2}, \gamma_{3}, t) \szp{]}   + \; \szp{[} F_{j2}( \gamma_{1}, \gamma_{2}, \gamma_{3}, t )\; -\;\frac{\gamma_{2}}{\gamma_{1}} \, F_{j1}( \gamma_{1}, \gamma_{2}, \gamma_{3}, t) \szp{]} 
\quad\quad\quad\quad\quad\quad 
\end{eqnarray}

\begin{eqnarray}\label{eqn142}
\frac{d}{dt} \szp{[}  U_{j3}( \gamma_{1}, \gamma_{2}, \gamma_{3}, t )\; -\;\frac{\gamma_{3}}{\gamma_{1}}\, U_{j1}( \gamma_{1}, \gamma_{2}, \gamma_{3}, t) \szp{]} \;= -\nu( \gamma_{1}^{2} +\gamma_{2}^{2} +\gamma_{3}^{2}) \szp{[} U_{j3}( \gamma_{1}, \gamma_{2}, \gamma_{3}, t )\; -\; \nonumber \\ 
\nonumber\\
-\;\frac{\gamma_{3}}{\gamma_{1}}\, U_{j1}( \gamma_{1}, \gamma_{2}, \gamma_{3}, t) \szp{]}   + \;\szp{[} F_{j3}( \gamma_{1}, \gamma_{2}, \gamma_{3}, t )\; -\;\frac{\gamma_{3}}{\gamma_{1}}\, F_{j1}( \gamma_{1}, \gamma_{2}, \gamma_{3}, t) \szp{]} 
\quad\quad\quad\quad\quad\quad 
\end{eqnarray}

\begin{equation}\label{eqn143}
\gamma_{1} U_{j1}( \gamma_{1}, \gamma_{2}, \gamma_{3}, t ) \;+\; \gamma_{2}\, U_{j2}( \gamma_{1}, \gamma_{2}, \gamma_{3}, t ) \;+\; \gamma_{3}\, U_{j3}( \gamma_{1}, \gamma_{2}, \gamma_{3}, t ) \;=\;0
\end{equation}

\begin{equation}\label{eqn144}
U_{j1}(\gamma_{1}, \gamma_{2},  \gamma_{3},  0)\;=\; U_{1}^{0}(\gamma_{1} ,\gamma_{2} ,\gamma_{3})
\end{equation}

\begin{equation}\label{eqn145}
U_{j2}(\gamma_{1}, \gamma_{2},  \gamma_{3},  0)\;=\; U_{2}^{0}(\gamma_{1} ,\gamma_{2} ,\gamma_{3})
\end{equation}

\begin{equation}\label{eqn146}
U_{j3}(\gamma_{1}, \gamma_{2},  \gamma_{3},  0)\;=\; U_{3}^{0}(\gamma_{1} ,\gamma_{2} ,\gamma_{3})
\end{equation}

On the third stage we can use Laplace transform $(\ref{A4}), (\ref{A5})$  for a system of linear ordinary differential equations with constant coefficients $(\ref{eqn141})\;-\; (\ref{eqn143})\;$  and have in result a system of linear algebraic equations with constant coefficients:

\[U_{jk}^{\otimes} (\gamma_{1}, \gamma_{2}, \gamma_{3}, \eta) \;=\;L[U_{jk}(\gamma_{1}, \gamma_{2}, \gamma_{3}, t)] \;\;\;\;\;\;\;      \rm{k=1,2,3}\] 

\begin{eqnarray}\label{eqn147}
\eta \szp{[}  U_{j2}^{\otimes}( \gamma_{1}, \gamma_{2}, \gamma_{3}, \eta )\; -\;\frac{\gamma_{2}}{\gamma_{1}} U_{j1}^{\otimes}( \gamma_{1}, \gamma_{2}, \gamma_{3}, \eta) \szp{]}  \;-\; \szp{[}   U_{j2}( \gamma_{1}, \gamma_{2}, \gamma_{3}, 0 )\; -\;\frac{\gamma_{2}}{\gamma_{1}} U_{j1}( \gamma_{1}, \gamma_{2}, \gamma_{3}, 0) \szp{]}  \;=  \nonumber
\\
\nonumber\\
-\nu( \gamma_{1}^{2} +\gamma_{2}^{2} +\gamma_{3}^{2})\szp{[}   U_{j2}^{\otimes}                          ( \gamma_{1}, \gamma_{2}, \gamma_{3}, \eta )\; -\;\frac{\gamma_{2}}{\gamma_{1}} U_{j1}^{\otimes}( \gamma_{1}, \gamma_{2}, \gamma_{3}, \eta) \szp{]}  \;+ 
\quad\quad\quad\quad\quad\quad 
\nonumber
\\
\nonumber\\
+\; \szp{[}   F_{j2}^{\otimes}( \gamma_{1}, \gamma_{2}, \gamma_{3}, \eta )\; -\;\frac{\gamma_{2}}{\gamma_{1}} F_{j1}^{\otimes}( \gamma_{1}, \gamma_{2}, \gamma_{3}, \eta) \szp{]}
\quad\quad\quad\quad\quad\quad \quad\quad\quad \quad\quad\quad 
\end{eqnarray}

\begin{eqnarray}\label{eqn148}
\eta  \szp{[}   U_{j3}^{\otimes}( \gamma_{1}, \gamma_{2}, \gamma_{3}, \eta )\; -\;\frac{\gamma_{3}}{\gamma_{1}} U_{j1}^{\otimes}( \gamma_{1}, \gamma_{2}, \gamma_{3}, \eta) \szp{]}   \;-\; \szp{[}    U_{j3}( \gamma_{1}, \gamma_{2}, \gamma_{3}, 0 )\; -\;\frac{\gamma_{3}}{\gamma_{1}} U_{j1}( \gamma_{1}, \gamma_{2}, \gamma_{3}, 0) \szp{]}   \;=  \nonumber
\\
\nonumber\\
-\nu( \gamma_{1}^{2} +\gamma_{2}^{2} +\gamma_{3}^{2})\szp{[}    U_{j3}^{\otimes}                          ( \gamma_{1}, \gamma_{2}, \gamma_{3}, \eta )\; -\;\frac{\gamma_{3}}{\gamma_{1}} U_{j1}^{\otimes}( \gamma_{1}, \gamma_{2}, \gamma_{3}, \eta) \szp{]}   \;+ 
\quad\quad\quad\quad\quad\quad 
\nonumber
\\
\nonumber\\
+\; \szp{[}    F_{j3}^{\otimes}( \gamma_{1}, \gamma_{2}, \gamma_{3}, \eta )\; -\;\frac{\gamma_{3}}{\gamma_{1}} F_{j1}^{\otimes}( \gamma_{1}, \gamma_{2}, \gamma_{3}, \eta) \szp{]} 
\quad\quad\quad\quad\quad\quad \quad\quad\quad \quad\quad\quad 
\end{eqnarray}

\begin{equation}\label{eqn149}
\gamma_{1} U_{j1}^{\otimes}( \gamma_{1}, \gamma_{2}, \gamma_{3}, \eta ) \;+\; \gamma_{2}\, U_{j2}^{\otimes}( \gamma_{1}, \gamma_{2}, \gamma_{3}, \eta ) \;+\; \gamma_{3}\, U_{j3}^{\otimes}( \gamma_{1}, \gamma_{2}, \gamma_{3}, \eta ) \;=\;0
\end{equation}

\begin{equation}\label{eqn150}
U_{j1}(\gamma_{1}, \gamma_{2},  \gamma_{3},  0)\;=\; U_{1}^{0}(\gamma_{1} ,\gamma_{2} ,\gamma_{3})
\end{equation}

\begin{equation}\label{eqn151}
U_{j2}(\gamma_{1}, \gamma_{2},  \gamma_{3},  0)\;=\; U_{2}^{0}(\gamma_{1} ,\gamma_{2} ,\gamma_{3})
\end{equation}

\begin{equation}\label{eqn152}
U_{j3}(\gamma_{1}, \gamma_{2},  \gamma_{3},  0)\;=\; U_{3}^{0}(\gamma_{1} ,\gamma_{2} ,\gamma_{3})
\end{equation}

Let us rewrite system of equations $(\ref{eqn147})\;-\; (\ref{eqn149})\;$ in the following form:

\begin{eqnarray}\label{eqn147a}
\szp{[}\eta \;+\;\nu( \gamma_{1}^{2} +\gamma_{2}^{2} +\gamma_{3}^{2})\szp{]}\frac{\gamma_{2}}{\gamma_{1}} U_{j1}^{\otimes}( \gamma_{1}, \gamma_{2}, \gamma_{3}, \eta)   \;-\; \szp{[}\eta \;+\;\nu( \gamma_{1}^{2} +\gamma_{2}^{2} +\gamma_{3}^{2})\szp{]} U_{j2}^{\otimes}( \gamma_{1}, \gamma_{2}, \gamma_{3}, \eta)  \;=  \nonumber
\\
\nonumber\\
\;\;\; \szp{[}  \frac{\gamma_{2}}{\gamma_{1}} F_{j1}^{\otimes}( \gamma_{1}, \gamma_{2}, \gamma_{3}, \eta) \; -\;F_{j2}^{\otimes}( \gamma_{1}, \gamma_{2}, \gamma_{3}, \eta ) \szp{]}+ \szp{[} \frac{\gamma_{2}}{\gamma_{1}} U_{j1}( \gamma_{1}, \gamma_{2}, \gamma_{3}, 0)  \; -\;U_{j2}( \gamma_{1}, \gamma_{2}, \gamma_{3}, 0 ) \szp{]}
\quad
\end{eqnarray}

\begin{eqnarray}\label{eqn148a}
\szp{[}\eta \;+\;\nu( \gamma_{1}^{2} +\gamma_{2}^{2} +\gamma_{3}^{2})\szp{]}\frac{\gamma_{3}}{\gamma_{1}} U_{j1}^{\otimes}( \gamma_{1}, \gamma_{2}, \gamma_{3}, \eta)   \;-\; \szp{[}\eta \;+\;\nu( \gamma_{1}^{2} +\gamma_{2}^{2} +\gamma_{3}^{2})\szp{]} U_{j3}^{\otimes}( \gamma_{1}, \gamma_{2}, \gamma_{3}, \eta)  \;=  \nonumber
\\
\nonumber\\
\;\;\; \szp{[}  \frac{\gamma_{3}}{\gamma_{1}} F_{j1}^{\otimes}( \gamma_{1}, \gamma_{2}, \gamma_{3}, \eta) \; -\;F_{j3}^{\otimes}( \gamma_{1}, \gamma_{2}, \gamma_{3}, \eta ) \szp{]}+ \szp{[} \frac{\gamma_{3}}{\gamma_{1}} U_{j1}( \gamma_{1}, \gamma_{2}, \gamma_{3}, 0)  \; -\;U_{j3}( \gamma_{1}, \gamma_{2}, \gamma_{3}, 0 ) \szp{]}
\quad
\end{eqnarray}

\begin{equation}\label{eqn149a}
\gamma_{1} U_{j1}^{\otimes}( \gamma_{1}, \gamma_{2}, \gamma_{3}, \eta ) \;+\; \gamma_{2}\, U_{j2}^{\otimes}( \gamma_{1}, \gamma_{2}, \gamma_{3}, \eta ) \;+\; \gamma_{3}\, U_{j3}^{\otimes}( \gamma_{1}, \gamma_{2}, \gamma_{3}, \eta ) \;=\;0
\end{equation}

Determinant of this system is

\begin{eqnarray}\label{eqn1481b}
\Delta = \left| \begin{array}{ccc}
\szp{[}\eta \;+\;\nu( \gamma_{1}^{2} +\gamma_{2}^{2} +\gamma_{3}^{2})\szp{]}\frac{\gamma_{2}}{\gamma_{1}} & -\szp{[}\eta \;+\;\nu( \gamma_{1}^{2} +\gamma_{2}^{2} +\gamma_{3}^{2})\szp{]} & 0 \\
\nonumber\\
 \szp{[}\eta \;+\;\nu( \gamma_{1}^{2} +\gamma_{2}^{2} +\gamma_{3}^{2})\szp{]}\frac{\gamma_{3}}{\gamma_{1}} & 0 & -\szp{[}\eta \;+\;\nu( \gamma_{1}^{2} +\gamma_{2}^{2} +\gamma_{3}^{2})\szp{]} \\
\nonumber\\ 
\gamma_{1} & \gamma_{2} & \gamma_{3} \end{array} \right| = 
\nonumber\\ 
\nonumber\\
\nonumber\\
=\;\frac{\szp{[}\eta \;+\;\nu( \gamma_{1}^{2} +\gamma_{2}^{2} +\gamma_{3}^{2})\szp{]}^{2}( \gamma_{1}^{2} +\gamma_{2}^{2} +\gamma_{3}^{2})}{\gamma_{1}}
\neq \;\; 0
\quad\quad\quad\quad\quad\quad\quad\quad\quad\quad\quad\quad
\end{eqnarray}

And consequently the system of equations $(\ref{eqn147})\;-\; (\ref{eqn149})\;$ and/or $\;(\ref{eqn147a})\;-\; (\ref{eqn149a})\;$ has a unique solution. Taking into account formulas $(\ref{eqn150})\;-\; (\ref{eqn152})\;$ we can write this solution in the following form:

\begin{eqnarray}\label{eqn153}
U_{j1}^{\otimes}( \gamma_{1}, \gamma_{2}, \gamma_{3}, \eta )\;=\;\frac{[( \gamma_{2}^{2} +\gamma_{3}^{2})  F_{j1}^{\otimes}( \gamma_{1}, \gamma_{2}, \gamma_{3}, \eta) - \gamma_{1}\gamma_{2} F_{j2}^{\otimes}( \gamma_{1}, \gamma_{2}, \gamma_{3}, \eta) - \gamma_{1}\gamma_{3} F_{j3}^{\otimes}( \gamma_{1}, \gamma_{2}, \gamma_{3}, \eta)]}{ (\gamma_{1}^{2} +\gamma_{2}^{2} +\gamma_{3}^{2}) [\eta+\nu (\gamma_{1}^{2} +\gamma_{2}^{2} +\gamma_{3}^{2})] }\;+\nonumber
\\
\nonumber\\
+\; \frac{ U_{1}^{0}(\gamma_{1} , \gamma_{2} , \gamma_{3})}{[\eta+\nu (\gamma_{1}^{2} +\gamma_{2}^{2} +\gamma_{3}^{2})] } 
\quad\quad\quad\quad\quad\quad \quad\quad\quad \quad\quad\quad 
\quad\quad\quad\quad\quad\quad 
\end{eqnarray}

\begin{eqnarray}\label{eqn154}
U_{j2}^{\otimes}( \gamma_{1}, \gamma_{2}, \gamma_{3}, \eta )\;=\;\frac{[( \gamma_{3}^{2} +\gamma_{1}^{2})  F_{j2}^{\otimes}( \gamma_{1}, \gamma_{2}, \gamma_{3}, \eta) - \gamma_{2}\gamma_{3} F_{j3}^{\otimes}( \gamma_{1}, \gamma_{2}, \gamma_{3}, \eta) - \gamma_{2}\gamma_{1} F_{j1}^{\otimes}( \gamma_{1}, \gamma_{2}, \gamma_{3}, \eta)]}{ (\gamma_{1}^{2} +\gamma_{2}^{2} +\gamma_{3}^{2}) [\eta+\nu (\gamma_{1}^{2} +\gamma_{2}^{2} +\gamma_{3}^{2})] }\;+\nonumber
\\
\nonumber\\
+\; \frac{ U_{2}^{0}(\gamma_{1} , \gamma_{2} , \gamma_{3})}{[\eta+\nu (\gamma_{1}^{2} +\gamma_{2}^{2} +\gamma_{3}^{2})] } 
\quad\quad\quad\quad\quad\quad \quad\quad\quad \quad\quad\quad 
\quad\quad\quad\quad\quad\quad 
\end{eqnarray}

\begin{eqnarray}\label{eqn155}
U_{j3}^{\otimes}( \gamma_{1}, \gamma_{2}, \gamma_{3}, \eta )\;=\;\frac{[( \gamma_{1}^{2} +\gamma_{2}^{2})  F_{j3}^{\otimes}( \gamma_{1}, \gamma_{2}, \gamma_{3}, \eta) - \gamma_{3}\gamma_{1} F_{j1}^{\otimes}( \gamma_{1}, \gamma_{2}, \gamma_{3}, \eta) - \gamma_{3}\gamma_{2} F_{j2}^{\otimes}( \gamma_{1}, \gamma_{2}, \gamma_{3}, \eta)]}{ (\gamma_{1}^{2} +\gamma_{2}^{2} +\gamma_{3}^{2}) [\eta+\nu (\gamma_{1}^{2} +\gamma_{2}^{2} +\gamma_{3}^{2})] }\;+ \nonumber
\\
\nonumber\\
+\; \frac{ U_{3}^{0}(\gamma_{1} , \gamma_{2} , \gamma_{3})}{[\eta+\nu (\gamma_{1}^{2} +\gamma_{2}^{2} +\gamma_{3}^{2})] } 
\quad\quad\quad\quad\quad\quad \quad\quad\quad \quad\quad\quad 
\quad\quad\quad\quad\quad\quad 
\end{eqnarray}

Then we use the convolution theorem with the convolution formula (\ref{A6}) and integral (\ref{A7}) for $(\ref{eqn153})\;-\; (\ref{eqn155})\;$and obtain:

\begin{eqnarray}\label{eqn156}
U_{j1}(\gamma_{1}, \gamma_{2}, \gamma_{3}, t)\;=\;
\quad\quad\quad\quad\quad\quad \quad\quad\quad \quad\quad\quad 
\quad\quad\quad\quad\quad\quad
 \nonumber
\\
\nonumber\\
\int_{0}^{t} \sz{e} ^{-\nu (\gamma_{1}^{2} +\gamma_{2}^{2} +\gamma_{3}^{2}) (t-\tau)} \frac{[( \gamma_{2}^{2} +\gamma_{3}^{2})  F_{j1}( \gamma_{1}, \gamma_{2}, \gamma_{3}, \tau) -\gamma_{1}\gamma_{2} F_{j2}( \gamma_{1}, \gamma_{2}, \gamma_{3}, \tau ) -\gamma_{1}\gamma_{3} F_{j3}( \gamma_{1}, \gamma_{2}, \gamma_{3}, \tau ) ]}{ (\gamma_{1}^{2} +\gamma_{2}^{2}+\gamma_{3}^{2} ) }\,d\tau\;+ \nonumber
\\
\nonumber\\
+\;\sz{e} ^{-\nu (\gamma_{1}^{2} +\gamma_{2}^{2} +\gamma_{3}^{2}) t} \;U_{1}^{0}(\gamma_{1} ,\gamma_{2} ,\gamma_{3})
\quad\quad\quad\quad\quad\quad \quad\quad\quad \quad\quad\quad 
\quad\quad\quad\quad 
\end{eqnarray}

\begin{eqnarray}\label{eqn157}
U_{j2}(\gamma_{1}, \gamma_{2}, \gamma_{3}, t)\;=\; 
\quad\quad\quad\quad\quad\quad \quad\quad\quad \quad\quad\quad 
\quad\quad\quad\quad\quad\quad
\nonumber
\\
\nonumber\\
\int_{0}^{t} \sz{e} ^{-\nu (\gamma_{1}^{2} +\gamma_{2}^{2} +\gamma_{3}^{2}) (t-\tau)} \frac{[( \gamma_{3}^{2} +\gamma_{1}^{2})  F_{j2}( \gamma_{1}, \gamma_{2}, \gamma_{3}, \tau) -\gamma_{2}\gamma_{3} F_{j3}( \gamma_{1}, \gamma_{2}, \gamma_{3}, \tau ) -\gamma_{2}\gamma_{1} F_{j1}( \gamma_{1}, \gamma_{2}, \gamma_{3}, \tau ) ]}{ (\gamma_{1}^{2} +\gamma_{2}^{2}+\gamma_{3}^{2} ) }\,d\tau\;+ \nonumber
\\
\nonumber\\
+\;\sz{e} ^{-\nu (\gamma_{1}^{2} +\gamma_{2}^{2} +\gamma_{3}^{2}) t} \;U_{2}^{0}(\gamma_{1} ,\gamma_{2} ,\gamma_{3})
\quad\quad\quad\quad\quad\quad \quad\quad\quad \quad\quad\quad 
\quad\quad\quad\quad 
\end{eqnarray}

\begin{eqnarray}\label{eqn158}
U_{j3}(\gamma_{1}, \gamma_{2}, \gamma_{3}, t)\;=\; 
\quad\quad\quad\quad\quad\quad \quad\quad\quad \quad\quad\quad 
\quad\quad\quad\quad\quad\quad
\nonumber
\\
\nonumber\\
\int_{0}^{t} \sz{e} ^{-\nu (\gamma_{1}^{2} +\gamma_{2}^{2} +\gamma_{3}^{2}) (t-\tau)} \frac{[( \gamma_{1}^{2} +\gamma_{2}^{2})  F_{j3}( \gamma_{1}, \gamma_{2}, \gamma_{3}, \tau) -\gamma_{3}\gamma_{1} F_{j1}( \gamma_{1}, \gamma_{2}, \gamma_{3}, \tau ) -\gamma_{3}\gamma_{2} F_{j2}( \gamma_{1}, \gamma_{2}, \gamma_{3}, \tau ) ]}{ (\gamma_{1}^{2} +\gamma_{2}^{2}+\gamma_{3}^{2} ) }\,d\tau\;+ 
\nonumber
\\
\nonumber\\
+\;\sz{e} ^{-\nu (\gamma_{1}^{2} +\gamma_{2}^{2} +\gamma_{3}^{2}) t} \;U_{3}^{0}(\gamma_{1} ,\gamma_{2} ,\gamma_{3})
\quad\quad\quad\quad\quad\quad \quad\quad\quad \quad\quad\quad 
\quad\quad\quad\quad 
\end{eqnarray}

$ P_{j}(\gamma_{1}, \gamma_{2}, \gamma_{3}, t) $ is obtained from equations $(\ref{eqn134}),\;(\ref{eqn135})\;$,$\; (\ref{eqn136})\;$ with use of equation $\; (\ref{eqn137})\;$:

\begin{equation}\label{eqn159}
P_{j}(\gamma_{1}, \gamma_{2}, \gamma_{3}, t) \;=\;i\frac{[\gamma_{1} F_{j1}( \gamma_{1}, \gamma_{2}, \gamma_{3}, t) +\gamma_{2} F_{j2}( \gamma_{1}, \gamma_{2}, \gamma_{3}, t ) +\gamma_{3} F_{j3}( \gamma_{1}, \gamma_{2}, \gamma_{3}, t )]}{ (\gamma_{1}^{2} +\gamma_{2}^{2} +\gamma_{3}^{2}) }
\end{equation}
\\
Using of the Fourier inversion formula $(\ref{A3})$ leads to:

\begin{eqnarray}\label{eqn160}
u_{j1}(x_{1}, x_{2}, x_{3}, t)\;=\; 
\frac{1}{(2\pi)^{3/2}} \int_{-\infty}^{\infty} \int_{-\infty}^{\infty} \int_{-\infty}^{\infty} \biggl[ \int_{0}^{t} \sz{e} ^{-\nu (\gamma_{1}^{2} +\gamma_{2}^{2} +\gamma_{3}^{2}) (t-\tau)} \frac{[ ( \gamma_{2}^{2} +\gamma_{3}^{2})  F_{j1}( \gamma_{1}, \gamma_{2}, \gamma_{3}, \tau)]} { (\gamma_{1}^{2} +\gamma_{2}^{2}+\gamma_{3}^{2} ) } \,d\tau \;- \nonumber
\\
\nonumber\\
-\; \int_{0}^{t} \sz{e} ^{-\nu (\gamma_{1}^{2} +\gamma_{2}^{2} +\gamma_{3}^{2}) (t-\tau)}\frac{ [\gamma_{1}\gamma_{2} F_{j2}( \gamma_{1}, \gamma_{2}, \gamma_{3}, \tau ) + \gamma_{1}\gamma_{3} F_{j3}( \gamma_{1}, \gamma_{2}, \gamma_{3}, \tau )]} { (\gamma_{1}^{2} +\gamma_{2}^{2}+\gamma_{3}^{2} ) }\,d\tau\;   +
\quad
\nonumber\\
\nonumber\\
+\;\sz{e} ^{-\nu (\gamma_{1}^{2} +\gamma_{2}^{2} +\gamma_{3}^{2}) t} \;U_{1}^{0}(\gamma_{1} ,\gamma_{2} ,\gamma_{3})\biggr] \;\sz{e} ^{-i(x_{1}\gamma_{1}+x_{2}\gamma_{2}+x_{3}\gamma_{3})}\,d\gamma_{1}d\gamma_{2}d\gamma_{3}\;=
\nonumber\\
\nonumber\\
=\;\frac{1}{8\pi^{3}} \int_{-\infty}^{\infty} \int_{-\infty}^{\infty} \int_{-\infty}^{\infty}\frac{( \gamma_{2}^{2} +\gamma_{3}^{2})} { (\gamma_{1}^{2} +\gamma_{2}^{2}+\gamma_{3}^{2} ) }  \biggl[ \int_{0}^{t} \sz{e} ^{-\nu (\gamma_{1}^{2} +\gamma_{2}^{2} +\gamma_{3}^{2}) (t-\tau)} \int_{-\infty}^{\infty}\int_{-\infty}^{\infty}\int_{-\infty}^{\infty}\sz{e}  ^{i(\tilde x_{1}\gamma_{1}+\tilde x_{2}\gamma_{2}+\tilde x_{3}\gamma_{3})} \cdot 
\nonumber\\
\nonumber\\
\cdot f_{j1}(\tilde x_{1},\tilde x_{2}, \tilde x_{3},\tau)\,d\tilde x_{1}d\tilde x_{2} d\tilde x_{3}d\tau\biggr]\sz{e} ^{-i(x_{1}\gamma_{1}+x_{2}\gamma_{2}+x_{3}\gamma_{3})}\,d\gamma_{1}d\gamma_{2}d\gamma_{3}\;-
\nonumber\\
\nonumber\\
-\;\frac{1}{8\pi^{3}} \int_{-\infty}^{\infty} \int_{-\infty}^{\infty} \int_{-\infty}^{\infty}\frac{ \gamma_{1}\gamma_{2}} { (\gamma_{1}^{2} +\gamma_{2}^{2}+\gamma_{3}^{2} ) }  \biggl[ \int_{0}^{t} \sz{e} ^{-\nu (\gamma_{1}^{2} +\gamma_{2}^{2} +\gamma_{3}^{2}) (t-\tau)} \int_{-\infty}^{\infty}\int_{-\infty}^{\infty}\int_{-\infty}^{\infty}\sz{e}  ^{i(\tilde x_{1}\gamma_{1}+\tilde x_{2}\gamma_{2}+\tilde x_{3}\gamma_{3})} \cdot 
\nonumber\\
\nonumber\\
\cdot f_{j2}(\tilde x_{1},\tilde x_{2}, \tilde x_{3},\tau)\,d\tilde x_{1}d\tilde x_{2} d\tilde x_{3}d\tau\biggr]\sz{e} ^{-i(x_{1}\gamma_{1}+x_{2}\gamma_{2}+x_{3}\gamma_{3})}\,d\gamma_{1}d\gamma_{2}d\gamma_{3}\;-
\nonumber\\
\nonumber\\
-\;\frac{1}{8\pi^{3}} \int_{-\infty}^{\infty} \int_{-\infty}^{\infty} \int_{-\infty}^{\infty}\frac{ \gamma_{1}\gamma_{3}} { (\gamma_{1}^{2} +\gamma_{2}^{2}+\gamma_{3}^{2} ) }  \biggl[ \int_{0}^{t} \sz{e} ^{-\nu (\gamma_{1}^{2} +\gamma_{2}^{2} +\gamma_{3}^{2}) (t-\tau)} \int_{-\infty}^{\infty}\int_{-\infty}^{\infty}\int_{-\infty}^{\infty}\sz{e}  ^{i(\tilde x_{1}\gamma_{1}+\tilde x_{2}\gamma_{2}+\tilde x_{3}\gamma_{3})} \cdot 
\nonumber\\
\nonumber\\
\cdot f_{j3}(\tilde x_{1},\tilde x_{2}, \tilde x_{3},\tau)\,d\tilde x_{1}d\tilde x_{2} d\tilde x_{3}d\tau\biggr]\sz{e} ^{-i(x_{1}\gamma_{1}+x_{2}\gamma_{2}+x_{3}\gamma_{3})}\,d\gamma_{1}d\gamma_{2}d\gamma_{3}\;+
\nonumber\\
\nonumber\\
+\;\frac{1}{8\pi^{3}} \int_{-\infty}^{\infty} \int_{-\infty}^{\infty} \int_{-\infty}^{\infty} \sz{e} ^{-\nu (\gamma_{1}^{2} +\gamma_{2}^{2} +\gamma_{3}^{2}) t}\biggl[  \int_{-\infty}^{\infty}\int_{-\infty}^{\infty}\int_{-\infty}^{\infty}\sz{e}  ^{i(\tilde x_{1}\gamma_{1}+\tilde x_{2}\gamma_{2}+\tilde x_{3}\gamma_{3})}
 \cdot \quad\quad\quad\quad\quad\quad\quad\quad\quad
\nonumber\\
\nonumber\\
\cdot \; u_{1}^{0}(\tilde x_{1},\tilde x_{2}, \tilde x_{3})\,d\tilde x_{1}d\tilde x_{2} d\tilde x_{3}\biggr]\sz{e} ^{-i(x_{1}\gamma_{1}+x_{2}\gamma_{2}+x_{3}\gamma_{3})}\,d\gamma_{1}d\gamma_{2}d\gamma_{3} \;= 
\nonumber\\
\nonumber\\
=\; S_{11}(f_{j1})\;+\; S_{12}(f_{j2})\;+\; S_{13}(f_{j3})\;+\;B(u_{1}^0) \quad\quad\quad\quad\quad\quad \quad\quad\quad \quad\quad\quad 
\end{eqnarray}

\begin{eqnarray}\label{eqn161}
u_{j2}(x_{1}, x_{2}, x_{3}, t)\;=\;
\frac{1}{(2\pi)^{3/2}} \int_{-\infty}^{\infty} \int_{-\infty}^{\infty} \int_{-\infty}^{\infty} \biggl[ \int_{0}^{t} \sz{e} ^{-\nu (\gamma_{1}^{2} +\gamma_{2}^{2} +\gamma_{3}^{2}) (t-\tau)} \frac{[( \gamma_{3}^{2} +\gamma_{1}^{2})  F_{j2}( \gamma_{1}, \gamma_{2}, \gamma_{3}, \tau)]} { (\gamma_{1}^{2} +\gamma_{2}^{2}+\gamma_{3}^{2} ) } \,d\tau \;-
\nonumber\\
\nonumber\\
-\; \int_{0}^{t} \sz{e} ^{-\nu (\gamma_{1}^{2} +\gamma_{2}^{2} +\gamma_{3}^{2}) (t-\tau)}\frac{[\gamma_{2}\gamma_{3} F_{j3}( \gamma_{1}, \gamma_{2}, \gamma_{3}, \tau ) +\gamma_{2}\gamma_{1} F_{j1}( \gamma_{1}, \gamma_{2}, \gamma_{3}, \tau )]} { (\gamma_{1}^{2} +\gamma_{2}^{2}+\gamma_{3}^{2} ) }\,d\tau\;   +
\quad
\nonumber\\
\nonumber\\
+\;\sz{e} ^{-\nu (\gamma_{1}^{2} +\gamma_{2}^{2} +\gamma_{3}^{2}) t} \;U_{2}^{0}(\gamma_{1} ,\gamma_{2} ,\gamma_{3})\biggr] \;\sz{e} ^{-i(x_{1}\gamma_{1}+x_{2}\gamma_{2}+x_{3}\gamma_{3})}\,d\gamma_{1}d\gamma_{2}d\gamma_{3}\;=
\nonumber\\
\nonumber\\
=\;-\;\frac{1}{8\pi^{3}} \int_{-\infty}^{\infty} \int_{-\infty}^{\infty} \int_{-\infty}^{\infty}\frac{ \gamma_{2}\gamma_{1}} { (\gamma_{1}^{2} +\gamma_{2}^{2}+\gamma_{3}^{2} ) }  \biggl[ \int_{0}^{t} \sz{e} ^{-\nu (\gamma_{1}^{2} +\gamma_{2}^{2} +\gamma_{3}^{2}) (t-\tau)} \int_{-\infty}^{\infty}\int_{-\infty}^{\infty}\int_{-\infty}^{\infty}\sz{e}  ^{i(\tilde x_{1}\gamma_{1}+\tilde x_{2}\gamma_{2}+\tilde x_{3}\gamma_{3})} \cdot 
\nonumber\\
\nonumber\\
\cdot f_{j1}(\tilde x_{1},\tilde x_{2}, \tilde x_{3},\tau)\,d\tilde x_{1}d\tilde x_{2} d\tilde x_{3}d\tau\biggr]\sz{e} ^{-i(x_{1}\gamma_{1}+x_{2}\gamma_{2}+x_{3}\gamma_{3})}\,d\gamma_{1}d\gamma_{2}d\gamma_{3}\;+
\nonumber\\
\nonumber\\
+\;\frac{1}{8\pi^{3}} \int_{-\infty}^{\infty} \int_{-\infty}^{\infty} \int_{-\infty}^{\infty}\frac{ (\gamma_{3}^2 + \gamma_{1}^2)} { (\gamma_{1}^{2} +\gamma_{2}^{2}+\gamma_{3}^{2} ) }  \biggl[ \int_{0}^{t} \sz{e} ^{-\nu (\gamma_{1}^{2} +\gamma_{2}^{2} +\gamma_{3}^{2}) (t-\tau)} \int_{-\infty}^{\infty}\int_{-\infty}^{\infty}\int_{-\infty}^{\infty}\sz{e}  ^{i(\tilde x_{1}\gamma_{1}+\tilde x_{2}\gamma_{2}+\tilde x_{3}\gamma_{3})} \cdot 
\nonumber\\
\nonumber\\
\cdot f_{j2}(\tilde x_{1},\tilde x_{2}, \tilde x_{3},\tau)\,d\tilde x_{1}d\tilde x_{2} d\tilde x_{3}d\tau\biggr]\sz{e} ^{-i(x_{1}\gamma_{1}+x_{2}\gamma_{2}+x_{3}\gamma_{3})}\,d\gamma_{1}d\gamma_{2}d\gamma_{3}\;-
\nonumber\\
\nonumber\\
-\;\frac{1}{8\pi^{3}} \int_{-\infty}^{\infty} \int_{-\infty}^{\infty} \int_{-\infty}^{\infty}\frac{ \gamma_{2}\gamma_{3}} { (\gamma_{1}^{2} +\gamma_{2}^{2}+\gamma_{3}^{2} ) }  \biggl[ \int_{0}^{t} \sz{e} ^{-\nu (\gamma_{1}^{2} +\gamma_{2}^{2} +\gamma_{3}^{2}) (t-\tau)} \int_{-\infty}^{\infty}\int_{-\infty}^{\infty}\int_{-\infty}^{\infty}\sz{e}  ^{i(\tilde x_{1}\gamma_{1}+\tilde x_{2}\gamma_{2}+\tilde x_{3}\gamma_{3})} \cdot 
\nonumber\\
\nonumber\\
\cdot f_{j3}(\tilde x_{1},\tilde x_{2}, \tilde x_{3},\tau)\,d\tilde x_{1}d\tilde x_{2} d\tilde x_{3}d\tau\biggr]\sz{e} ^{-i(x_{1}\gamma_{1}+x_{2}\gamma_{2}+x_{3}\gamma_{3})}\,d\gamma_{1}d\gamma_{2}d\gamma_{3}\;+
\nonumber\\
\nonumber\\
+\;\frac{1}{8\pi^{3}} \int_{-\infty}^{\infty} \int_{-\infty}^{\infty} \int_{-\infty}^{\infty} \sz{e} ^{-\nu (\gamma_{1}^{2} +\gamma_{2}^{2} +\gamma_{3}^{2}) t}\biggl[  \int_{-\infty}^{\infty}\int_{-\infty}^{\infty}\int_{-\infty}^{\infty}\sz{e}  ^{i(\tilde x_{1}\gamma_{1}+\tilde x_{2}\gamma_{2}+\tilde x_{3}\gamma_{3})}
 \cdot \quad\quad\quad\quad\quad\quad\quad\quad\quad
\nonumber\\
\nonumber\\
\cdot \; u_{2}^{0}(\tilde x_{1},\tilde x_{2}, \tilde x_{3})\,d\tilde x_{1}d\tilde x_{2} d\tilde x_{3}\biggr]\sz{e} ^{-i(x_{1}\gamma_{1}+x_{2}\gamma_{2}+x_{3}\gamma_{3})}\,d\gamma_{1}d\gamma_{2}d\gamma_{3} \;= 
\nonumber\\
\nonumber\\
=\; S_{21}(f_{j1})\;+\; S_{22}(f_{j2})\;+\; S_{23}(f_{j3})\;+\;B(u_{2}^0) 
\quad\quad\quad\quad\quad\quad \quad\quad\quad \quad\quad\quad 
\end{eqnarray}

\begin{eqnarray}\label{eqn162}
u_{j3}(x_{1}, x_{2}, x_{3}, t)\;=\;
\frac{1}{(2\pi)^{3/2}} \int_{-\infty}^{\infty} \int_{-\infty}^{\infty} \int_{-\infty}^{\infty} \biggl[ \int_{0}^{t} \sz{e} ^{-\nu (\gamma_{1}^{2} +\gamma_{2}^{2} +\gamma_{3}^{2}) (t-\tau)} \frac{[( \gamma_{1}^{2} +\gamma_{2}^{2})  F_{j3}( \gamma_{1}, \gamma_{2}, \gamma_{3}, \tau)]} { (\gamma_{1}^{2} +\gamma_{2}^{2}+\gamma_{3}^{2} ) } \,d\tau \;-
\nonumber\\
\nonumber\\
-\;\int_{0}^{t} \sz{e} ^{-\nu (\gamma_{1}^{2} +\gamma_{2}^{2} +\gamma_{3}^{2}) (t-\tau)}\frac{[\gamma_{3}\gamma_{1} F_{j1}( \gamma_{1}, \gamma_{2}, \gamma_{3}, \tau ) +\gamma_{3}\gamma_{2} F_{j2}( \gamma_{1}, \gamma_{2}, \gamma_{3}, \tau )]} { (\gamma_{1}^{2} +\gamma_{2}^{2}+\gamma_{3}^{2} ) }\,d\tau\;   +
\quad
\nonumber\\
\nonumber\\
+\; \sz{e} ^{-\nu (\gamma_{1}^{2} +\gamma_{2}^{2} +\gamma_{3}^{2}) t} \;U_{3}^{0}(\gamma_{1} ,\gamma_{2} ,\gamma_{3})\biggr] \;\sz{e} ^{-i(x_{1}\gamma_{1}+x_{2}\gamma_{2}+x_{3}\gamma_{3})}\,d\gamma_{1}d\gamma_{2}d\gamma_{3}\;=
\nonumber\\
\nonumber\\
=\;-\;\frac{1}{8\pi^{3}} \int_{-\infty}^{\infty} \int_{-\infty}^{\infty} \int_{-\infty}^{\infty}\frac{ \gamma_{3}\gamma_{1}} { (\gamma_{1}^{2} +\gamma_{2}^{2}+\gamma_{3}^{2} ) }  \biggl[ \int_{0}^{t} \sz{e} ^{-\nu (\gamma_{1}^{2} +\gamma_{2}^{2} +\gamma_{3}^{2}) (t-\tau)} \int_{-\infty}^{\infty}\int_{-\infty}^{\infty}\int_{-\infty}^{\infty}\sz{e}  ^{i(\tilde x_{1}\gamma_{1}+\tilde x_{2}\gamma_{2}+\tilde x_{3}\gamma_{3})} \cdot 
\nonumber\\
\nonumber\\
\cdot f_{j1}(\tilde x_{1},\tilde x_{2}, \tilde x_{3},\tau)\,d\tilde x_{1}d\tilde x_{2} d\tilde x_{3}d\tau\biggr]\sz{e} ^{-i(x_{1}\gamma_{1}+x_{2}\gamma_{2}+x_{3}\gamma_{3})}\,d\gamma_{1}d\gamma_{2}d\gamma_{3}\;-
\nonumber\\
\nonumber\\
-\;\frac{1}{8\pi^{3}} \int_{-\infty}^{\infty} \int_{-\infty}^{\infty} \int_{-\infty}^{\infty}\frac{ \gamma_{3}\gamma_{2}} { (\gamma_{1}^{2} +\gamma_{2}^{2}+\gamma_{3}^{2} ) }  \biggl[ \int_{0}^{t} \sz{e} ^{-\nu (\gamma_{1}^{2} +\gamma_{2}^{2} +\gamma_{3}^{2}) (t-\tau)} \int_{-\infty}^{\infty}\int_{-\infty}^{\infty}\int_{-\infty}^{\infty}\sz{e}  ^{i(\tilde x_{1}\gamma_{1}+\tilde x_{2}\gamma_{2}+\tilde x_{3}\gamma_{3})} \cdot 
\nonumber\\
\nonumber\\
\cdot f_{j2}(\tilde x_{1},\tilde x_{2}, \tilde x_{3},\tau)\,d\tilde x_{1}d\tilde x_{2} d\tilde x_{3}d\tau\biggr]\sz{e} ^{-i(x_{1}\gamma_{1}+x_{2}\gamma_{2}+x_{3}\gamma_{3})}\,d\gamma_{1}d\gamma_{2}d\gamma_{3}\;+
\nonumber\\
\nonumber\\
+\;\frac{1}{8\pi^{3}} \int_{-\infty}^{\infty} \int_{-\infty}^{\infty} \int_{-\infty}^{\infty}\frac{ (\gamma_{1}^2 + \gamma_{2}^2)} { (\gamma_{1}^{2} +\gamma_{2}^{2}+\gamma_{3}^{2} ) }  \biggl[ \int_{0}^{t} \sz{e} ^{-\nu (\gamma_{1}^{2} +\gamma_{2}^{2} +\gamma_{3}^{2}) (t-\tau)} \int_{-\infty}^{\infty}\int_{-\infty}^{\infty}\int_{-\infty}^{\infty}\sz{e}  ^{i(\tilde x_{1}\gamma_{1}+\tilde x_{2}\gamma_{2}+\tilde x_{3}\gamma_{3})} \cdot 
\nonumber\\
\nonumber\\
\cdot f_{j3}(\tilde x_{1},\tilde x_{2}, \tilde x_{3},\tau)\,d\tilde x_{1}d\tilde x_{2} d\tilde x_{3}d\tau\biggr]\sz{e} ^{-i(x_{1}\gamma_{1}+x_{2}\gamma_{2}+x_{3}\gamma_{3})}\,d\gamma_{1}d\gamma_{2}d\gamma_{3}\;+
\nonumber\\
\nonumber\\
+\;\frac{1}{8\pi^{3}} \int_{-\infty}^{\infty} \int_{-\infty}^{\infty} \int_{-\infty}^{\infty} \sz{e} ^{-\nu (\gamma_{1}^{2} +\gamma_{2}^{2} +\gamma_{3}^{2}) t}\biggl[  \int_{-\infty}^{\infty}\int_{-\infty}^{\infty}\int_{-\infty}^{\infty}\sz{e}  ^{i(\tilde x_{1}\gamma_{1}+\tilde x_{2}\gamma_{2}+\tilde x_{3}\gamma_{3})}
 \cdot \quad\quad\quad\quad\quad\quad\quad\quad\quad
\nonumber\\
\nonumber\\
\cdot \; u_{3}^{0}(\tilde x_{1},\tilde x_{2}, \tilde x_{3})\,d\tilde x_{1}d\tilde x_{2} d\tilde x_{3}\biggr]\sz{e} ^{-i(x_{1}\gamma_{1}+x_{2}\gamma_{2}+x_{3}\gamma_{3})}\,d\gamma_{1}d\gamma_{2}d\gamma_{3} \;= 
\nonumber\\
\nonumber\\
=\; S_{31}(f_{j1})\;+\; S_{32}(f_{j2})\;+\; S_{33}(f_{j3})\;+\;B(u_{3}^0) 
\quad\quad\quad\quad\quad\quad \quad\quad\quad \quad\quad\quad 
\end{eqnarray}

\begin{eqnarray}\label{eqn163}
p_{j}\,(x_{1}, x_{2}, x_{3}, t)\;=\; 
\frac{i}{(2\pi)^{3/2}} \int_{-\infty}^{\infty} \int_{-\infty}^{\infty} \int_{-\infty}^{\infty}\biggl[ \;\frac{[\gamma_{1} F_{j1}( \gamma_{1}, \gamma_{2}, \gamma_{3}, t) + \gamma_{2} F_{j2}( \gamma_{1}, \gamma_{2}, \gamma_{3}, t )]} { (\gamma_{1}^{2} +\gamma_{2}^{2} +\gamma_{3}^{2}) }\;+ 
\nonumber\\
\nonumber\\
+\;\frac{ \gamma_{3} F_{j3}( \gamma_{1}, \gamma_{2}, \gamma_{3}, t )}{  (\gamma_{1}^{2} +\gamma_{2}^{2} +\gamma_{3}^{2}) }\; \biggr] \;\sz{e} ^{-i(x_{1}\gamma_{1} +x_{2}\gamma_{2} +x_{3}\gamma_{3})}\,d\gamma_{1} d\gamma_{2} d\gamma_{3}\; = \; 
\nonumber\\
\nonumber\\
=\frac{i}{8\pi^{3}} \int_{-\infty}^{\infty} \int_{-\infty}^{\infty} \int_{-\infty}^{\infty}\frac{ \gamma_{1}} { (\gamma_{1}^{2} +\gamma_{2}^{2}+\gamma_{3}^{2} ) }  \biggl[ \int_{-\infty}^{\infty}\int_{-\infty}^{\infty}\int_{-\infty}^{\infty}\sz{e}  ^{i(\tilde x_{1}\gamma_{1}+\tilde x_{2}\gamma_{2}+\tilde x_{3}\gamma_{3})} \cdot 
\nonumber\\
\nonumber\\
\cdot f_{j1}(\tilde x_{1},\tilde x_{2}, \tilde x_{3},t)\,d\tilde x_{1}d\tilde x_{2} d\tilde x_{3}\biggr]\sz{e} ^{-i(x_{1}\gamma_{1}+x_{2}\gamma_{2}+x_{3}\gamma_{3})}\,d\gamma_{1}d\gamma_{2}d\gamma_{3}\;+
\nonumber\\
\nonumber\\
+\;\frac{i}{8\pi^{3}} \int_{-\infty}^{\infty} \int_{-\infty}^{\infty} \int_{-\infty}^{\infty}\frac{ \gamma_{2}} { (\gamma_{1}^{2} +\gamma_{2}^{2}+\gamma_{3}^{2} ) } \biggl [ \int_{-\infty}^{\infty}\int_{-\infty}^{\infty}\int_{-\infty}^{\infty}\sz{e}  ^{i(\tilde x_{1}\gamma_{1}+\tilde x_{2}\gamma_{2}+\tilde x_{3}\gamma_{3})} \cdot 
\nonumber\\
\nonumber\\
\cdot f_{j2}(\tilde x_{1},\tilde x_{2}, \tilde x_{3},t)\,d\tilde x_{1}d\tilde x_{2} d\tilde x_{3}\biggr]\sz{e} ^{-i(x_{1}\gamma_{1}+x_{2}\gamma_{2}+x_{3}\gamma_{3})}\,d\gamma_{1}d\gamma_{2}d\gamma_{3}\;+
\nonumber\\
\nonumber\\
+\;\frac{i}{8\pi^{3}} \int_{-\infty}^{\infty} \int_{-\infty}^{\infty} \int_{-\infty}^{\infty}\frac{ \gamma_{3}} { (\gamma_{1}^{2} +\gamma_{2}^{2}+\gamma_{3}^{2} ) }  \biggl[ \int_{-\infty}^{\infty}\int_{-\infty}^{\infty}\int_{-\infty}^{\infty}\sz{e}  ^{i(\tilde x_{1}\gamma_{1}+\tilde x_{2}\gamma_{2}+\tilde x_{3}\gamma_{3})} \cdot 
\nonumber\\
\nonumber\\
\cdot f_{j3}(\tilde x_{1},\tilde x_{2}, \tilde x_{3},t)\,d\tilde x_{1}d\tilde x_{2} d\tilde x_{3}\biggr]\sz{e} ^{-i(x_{1}\gamma_{1}+x_{2}\gamma_{2}+x_{3}\gamma_{3})}\,d\gamma_{1}d\gamma_{2}d\gamma_{3}\;=
\nonumber\\
\nonumber\\
=\;\tilde S_{1}(f_{j1})\;+\; \tilde S_{2}(f_{j2})\;+\; \tilde S_{3}(f_{j3})
\quad\quad\quad\quad\quad\quad \quad\quad \quad\quad\quad\quad\quad\quad\quad\quad \quad\quad  
\end{eqnarray}

So, the integrals $(\ref{eqn160})\;-\; (\ref{eqn163})\;$ exist [look below in \textbf {a)}, \textbf {b)} - the Fourier transform of the class S, infinitely differentiable functions]. For these calculations inverse
Fourier transforms are defined as Cauchy principal values and for $\gamma_{1}$=0, $\gamma_{2}$=0,$\gamma_{3}$=0.

Here $S_{11}(), S_{12}(), S_{13}(), S_{21}(), S_{22}(), S_{23}(), S_{31}(), S_{32}(), S_{33}(), B(), \tilde S_{1}(), \tilde S_{2}(), \tilde S_{3}()$ are 

the integral - operators.

\[S_{12}()\;= \;S_{21}() \]\[S_{13}()\;= \;S_{31}() \]\[S_{23}()\;= \;S_{32}() \]

We have for the vector $\vec{u}_{j}$ from the equations $(\ref{eqn160})\;-\; (\ref{eqn162})\;$:

\begin{equation}\label{eqn164}
\vec{u}_{j}\;=\;\bar{\bar{S}}\;\cdot\;\vec{f}_{j}\;+\;\bar{\bar{B}}\cdot\vec{u}^{0}\;,
\end{equation}

where $\;\bar{\bar{S}} \; $ and $\bar{\bar{B}}$ are the matrix integral operators:
\[ \left( \begin{array}{ccc}
S_{11} & S_{12} & S_{13} \\
S_{21} & S_{22} & S_{23} \\
S_{31} & S_{32} & S_{33}  \end{array} \right)\]
\[ \left( \begin{array}{ccc}
B & 0 & 0 \\
0 & B & 0 \\
0 & 0 & B  \end{array} \right)\]

and for the funcion $p_{j}$ from the equation $(\ref{eqn163})\;$:

\begin{equation}\label{eqn164a}
p_{j}\;=\;\tilde{\tilde{S}}\;\cdot\;\vec{f}_{j}\;,
\end{equation}

where $\;\tilde{\tilde{S}} \; $ is the matrix integral operator:
\[ \left( \begin{array}{ccc}
\tilde S_{1} & 0 & 0 \\
0 & \tilde S_{2} & 0 \\
0 & 0 & \tilde S_{3}  \end{array} \right)\]

We put $\vec{f}_{j}$ from equation $(\ref{eqn20})$ into equation $(\ref{eqn164})$ and have:

\begin{eqnarray}\label{eqn165}
\vec{u}_{j} = \bar{\bar{S}}\cdot(\;\vec{f}\;-\;(\;\vec{u}_{j-1}\cdot\nabla)\vec{u}_{j-1})\;+\;\bar{\bar{B}}\cdot\vec{u}^{0}\;=
\nonumber\\
\nonumber\\
=\;\bar{\bar{S}}\cdot\vec{f} \;-\;\bar{\bar{S}}\cdot(\vec{u}_{j-1}\;\cdot\;\nabla\;)\;\vec{u}_{j-1}\;+\;\bar{\bar{B}}\cdot\vec{u}^{0}\; =
\nonumber\\
\nonumber\\
=\;\vec{u}_{1}\;-\;\bar{\bar{S}}\cdot(\vec{u}_{j-1}\;\cdot\;\nabla)\;\vec{u}_{j-1}
\quad\quad\quad\quad\quad\quad \quad\quad 
\end{eqnarray}

Here $\vec{u}_{1}\;$ is the solution of the system of equations $(\ref{eqn13})\; - \;(\ref{eqn20})$ with condition:

\[\sum_{n=1}^{3} u_{n}\frac{\partial u_{k}}{\partial x_{n}}\;=\;0\;\;\;\;\;\;\;       \rm{k=1,2,3}\;\;\] 

For j = 1 formula $(\ref{eqn164})$ can be written as follows:

\begin{equation}\label{eqn168}
\vec{u}_{1}\;=\;\bar{\bar{S}}\;\cdot\;\vec{f}_{1}\;+\;\bar{\bar{B}}\cdot\vec{u}^{0}\;,\;\;\;\;\;\\\vec{f}_{1}(x,t)\; = \; \vec{f}(x,t)
\end{equation}

If t $\rightarrow$ 0 then $\vec{u}_{1} \rightarrow \vec{u}^{0}$ (look at matrix integral operators $\bar{\bar{S}}, \bar{\bar{B}}()\;\;$- integrals $\;(\ref{eqn160})\; - \;(\ref{eqn162})$).

For j = 2  we define from equation $(\ref{eqn20})$:

\begin{equation}
\vec{f}_{2}(x,t)\; = \; \vec{f}_{1}(x,t) \; - \;(\;\vec{u}_{1}\;\cdot\;\nabla\;)\;\vec{u}_{1}\;
\end{equation}

We denote: 

\begin{equation}\label{eqn169}
{\vec{f}_{2}^{*}\;=\;(\vec{u}_{1}\;\cdot\;\nabla)\;\vec{u}_{1}}
\end{equation}

and then we have: 

\begin{equation}
\vec{f}_{2}(x,t)\; = \; \vec{f}_{1}(x,t) \; - \vec{f}_{2}^{*}
\end{equation}

Then we get $\vec{u}_{2}$ from $(\ref{eqn164}),(\ref{eqn168})$:

\begin{equation}\label{eqn171}
\vec{u}_{2}\;=\;\bar{\bar{S}}\;\cdot\;\vec{f}_{2}  \;+\;\bar{\bar{B}}\cdot\vec{u}^{0}\;=\;\bar{\bar{S}}\;\cdot\;(\vec{f}_{1}\;-\;\vec{f}_{2}^{*})  \;+\;\bar{\bar{B}}\cdot\vec{u}^{0}\;=\;\vec{u}_{1}\;-\;\vec{u}_{2}^{*}
\end{equation}

Here we have:

\begin{equation}\label{eqn170}
{\vec{u}_{2}^{*}\;=\;\bar{\bar{S}}\;\cdot\;\vec{f}_{2}^{*}}
\end{equation}

If t $\rightarrow$ 0 then $\vec{u}_{2}^{*} \rightarrow$ 0 (look at matrix integral operator $\bar{\bar{S}}\;\;$- integrals $\;(\ref{eqn160})\; - \;(\ref{eqn162})$).

Continue for j = 3. We define from equation $(\ref{eqn20})$:

\begin{equation}
\vec{f}_{3}(x,t)\; = \; \vec{f}_{1}(x,t) \; - \;(\;\vec{u}_{2}\;\cdot\;\nabla\;)\;\vec{u}_{2}\;
\end{equation}

Here we have:

\begin{equation}\label{eqn172}
(\vec{u}_{2}\;\cdot\;\nabla)\;\vec{u}_{2}\;=
\;((\vec{u}_{1}\;-\;\vec{u}_{2}^{*})\;\cdot\;\nabla\;)
\;(\vec{u}_{1}\;-\;\vec{u}_{2}^{*})\;=
\;\vec{f}_{2}^{*}\;+\;\vec{f}_{3}^{*}
\end{equation}

We denote in $(\ref{eqn172})$:

\begin{equation}\label{eqn172z}
{\vec{f}_{3}^{*}\;=\;-\;  (\vec{u}_{1}\;\cdot\;\nabla)\;\vec{u}_{2}^{*}\; -\; (\vec{u}_{2}^{*}\;\cdot\;\nabla)\;\vec{u}_{1}\;+\;
(\vec{u}_{2}^{*}\;\cdot\;\nabla)\;\vec{u}_{2}^{*}}
\end{equation}

and then we have: 

\begin{equation}
\vec{f}_{3}(x,t)\; = \; \vec{f}_{1}(x,t) \; - \vec{f}_{2}^{*}\; - \vec{f}_{3}^{*}
\end{equation}

Then we get $\vec{u}_{3}$ from $(\ref{eqn164}) , (\ref{eqn168}) ,(\ref{eqn170})$:

\begin{equation}\label{eqn174}
\vec{u}_{3}\;=\;\bar{\bar{S}}\;\cdot\;\vec{f}_{3} \;+\;\bar{\bar{B}}\cdot\vec{u}^{0}\;=\;\bar{\bar{S}}\;\cdot\;(\vec{f}_{1} \;-\;\vec{f}_{2}^{*}\;-\;\vec{f}_{3}^{*})\;+\;\bar{\bar{B}}\cdot\vec{u}^{0}\;=\;\vec{u}_{1} \;-\;\vec{u}_{2}^{*}\;-\;\vec{u}_{3}^{*}
\end{equation}

Here we denote:

\begin{equation}\label{eqn173}
{\vec{u}_{3}^{*}\;=\;\bar{\bar{S}}\;\cdot\;\vec{f}_{3}^{*}}
\end{equation}

If t $\rightarrow$ 0 then $\vec{u}_{3}^{*} \rightarrow$ 0 (look at matrix integral operator $\bar{\bar{S}}\;\;$- integrals $\;(\ref{eqn160})\; - \;(\ref{eqn162})$).

For j = 4. We define from equation $(\ref{eqn20})$:

\begin{equation}
\vec{f}_{4}(x,t)\; = \; \vec{f}_{1}(x,t) \; - \;(\;\vec{u}_{3}\;\cdot\;\nabla\;)\;\vec{u}_{3}\;
\end{equation}

Here we have:

\begin{equation}\label{eqn175}
(\vec{u}_{3}\;\cdot\;\nabla)\;\vec{u}_{3}\;=
\;((\vec{u}_{2}\;-\;\vec{u}_{3}^{*})\;\cdot\;\nabla\;)
\;(\vec{u}_{2}\;-\;\vec{u}_{3}^{*})\;=
\;\vec{f}_{2}^{*}\;+\;\vec{f}_{3}^{*}\;+\;\vec{f}_{4}^{*}
\end{equation}

We denote in $(\ref{eqn175})$:

\begin{equation}\label{eqn175z}
{\vec{f}_{4}^{*}\;=\;-\;  (\vec{u}_{2}\;\cdot\;\nabla)\;\vec{u}_{3}^{*}\; -\; (\vec{u}_{3}^{*}\;\cdot\;\nabla)\;\vec{u}_{2}\;+\;
(\vec{u}_{3}^{*}\;\cdot\;\nabla)\;\vec{u}_{3}^{*}}
\end{equation}

and then we have: 

\begin{equation}
\vec{f}_{4}(x,t)\; = \; \vec{f}_{1}(x,t) \; - \vec{f}_{2}^{*}\; - \vec{f}_{3}^{*}\; - \vec{f}_{4}^{*}
\end{equation}

Then we get $\vec{u}_{4}$ from $(\ref{eqn164}) , (\ref{eqn168}) ,(\ref{eqn170}) ,(\ref{eqn173})$:

\begin{equation}\label{eqn177}
\vec{u}_{4}\;=\;\bar{\bar{S}}\;\cdot\;(\vec{f}_{1} \; -\;\vec{f}_{2}^{*}\; -\;\vec{f}_{3}^{*}\;-\;\vec{f}_{4}^{*})\;+\;\bar{\bar{B}}\cdot\vec{u}^{0}\;=\;\vec{u}_{1} \;-\;\vec{u}_{2}^{*}\;-\;\vec{u}_{3}^{*}\;-\;\vec{u}_{4}^{*}
\end{equation}

Here we denote:

\begin{equation}\label{eqn176}
{\vec{u}_{4}^{*}\;=\;\bar{\bar{S}}\;\cdot\;\vec{f}_{4}^{*}}
\end{equation}

If t $\rightarrow$ 0 then $\vec{u}_{4}^{*} \rightarrow$ 0 (look at matrix integral operator $\bar{\bar{S}}\;\;$- integrals $\;(\ref{eqn160})\; - \;(\ref{eqn162})$).

For arbitrary number $j$ $(j \geq 2)$. We define from equation $(\ref{eqn20})$:

\begin{equation}
\vec{f}_{j}(x,t)\; = \; \vec{f}_{1}(x,t) \; - \;(\;\vec{u}_{j-1}\;\cdot\;\nabla\;)\;\vec{u}_{j-1}\;
\end{equation}

Here we have:

\begin{equation}\label{eqn181}
(\vec{u}_{j-1}\;\cdot\;\nabla)\;\vec{u}_{j-1}\;=\;\sum_{l=2}^{j} \vec{f}_{l}^{*}
\end{equation}

and it follows:

\begin{equation}\label{eqn181a}
\vec{f}_{j}\;=\;\vec{f}_{1}\;-\; \sum_{l=2}^{j} \vec{f}_{l}^{*}
\end{equation}

{Here $\vec{f}_{2}^{*}$ is taken from formula $(\ref{eqn169}) $ and }

\begin{equation}\label{eqn181az}
{\vec{f}_{l}^{*}\;=\;-\;  (\vec{u}_{l-2}\;\cdot\;\nabla)\;\vec{u}_{l-1}^{*}\; -\; (\vec{u}_{l-1}^{*}\;\cdot\;\nabla)\;\vec{u}_{l-2}\;+\;
(\vec{u}_{l-1}^{*}\;\cdot\;\nabla)\;\vec{u}_{l-1}^{*}\;\;\;\;\;\;\;\;\;\;\;\;(l > 2) }
\end{equation}

Then we get $\vec{u}_{j}$ from $(\ref{eqn164}) , (\ref{eqn168}) $

\begin{equation}\label{eqn182}
\vec{u}_{j}\;=\;\bar{\bar{S}}\;\cdot\; \vec{f}_{j}\;+ \;\bar{\bar{B}}\cdot\vec{u}^{0}\;=\;\bar{\bar{S}}\;\cdot\;( \vec{f}_{1}\;-\; \sum_{l=2}^{j} \vec{f}_{l}^{*})\;+ \;\bar{\bar{B}}\cdot\vec{u}^{0}\;=\;\vec{u}_{1}\;-\;\sum_{l=2}^{j} \vec{u}_{l}^{*}\
\end{equation}

Here we denote:
\begin{equation}\label{eqn183}
{\vec{u}_{l}^{*}\;=\;\bar{\bar{S}}\;\cdot\;\vec{f}_{l}^{*}\;\;\;\;\;\;\;\;\;\;\;\;\;\;(2\;\leq\;l\;\leq\;j)}
\end{equation}

If t $\rightarrow$ 0 then $\vec{u}_{l}^{*} \rightarrow$ 0 (look at matrix integral operator $\bar{\bar{S}}\;\;$- integrals $\;(\ref{eqn160})\; - \;(\ref{eqn162})$).

We consider the equations $(\ref{eqn168})$ - $(\ref{eqn183})$ and see that the series $(\ref{eqn182})$ converges for $j \rightarrow \infty$ 

with the conditions for the first step (j = 1) of the iterative process:

\[\;\;\;\sum_{n=1}^{3} u_{0n}\frac{\partial u_{0k}}{\partial x_{n}} = 0\;\;\;\;\;\;\;\;       \rm{k=1,2,3}\]

Hence, we receive from equation $(\ref{eqn165})\;$ when $j \rightarrow \infty$:

\begin{equation}\label{eqn185}
\vec{u}_{\infty}=\;\vec{u}_{1}\;-\;\bar{\bar{S}}\cdot(\vec{u}_{\infty}\;\cdot\;\nabla)\;\vec{u}_{\infty}
\end{equation}

Equation $(\ref{eqn185})$ describes the converging iterative process. 
\nonumber\\

Below we show a proof that the iterative process is converging.

\textbf {a)} Let us consider S to be the class of all infinitely differentiable functions $\varphi(x)\; (-\infty < x < \infty)$,

 satisfying inequalities of the form 
 
\begin{equation}\label{eqn5a}
\mid x^{k}\varphi^{(q)}(x)\mid\;\leq\;C_{kq}\rm{\; for\; any }\;k, q = 0, 1, 2, ...
\end{equation}

where $C_{kq}$ is a constant and depends on $\varphi(x)$.

Then, \textbf {F}S = S, i.e., the Fourier transform operator \textbf {F} maps the class S onto the whole class S $\cite{gS01}$.

Now let us rewrite conditions $(\ref{eqn17})\;, \;(\ref{eqn18})$ in the following form:

\begin{equation}\label{eqn17a}
\mid(1+\mid \tilde x \mid)^{K}\partial_{\tilde x}^{\alpha}\vec{u}^{0}(\tilde x)\mid\;\leq\;C_{\alpha K} \quad\rm{on }\;R^{N}\;\rm{ for \; any }\;\alpha\;,\;K
\end{equation}

\begin{equation}\label{eqn18a}
\mid(1+\mid \tilde x \mid +\tau)^{K}\partial_{\tilde x}^{\alpha}\partial_{\tau}^{\beta}\vec{f}(\tilde x,\tau)\mid\;\leq\;C_{\alpha \beta K} \quad\rm{on }\;R^{N}\times[0,\infty)\;\rm{ for \; any }\;\alpha\;,\;\beta\;,\;K
\end{equation}

or for arbitrary k (${1\leq k \leq N}$)

\begin{equation}\label{eqn17b}
\mid(1+\mid \tilde x \mid)^{K}\partial_{\tilde x}^{\alpha}u_{k}^{0}(\tilde x)\mid\;\leq\;C_{\alpha K} \quad\rm{on }\;R^{N}\;\rm{ for \; any }\;\alpha\;,\;K
\end{equation}

\begin{equation}\label{eqn18b}
\mid(1+\mid \tilde x \mid +\tau)^{K}\partial_{\tilde x}^{\alpha}\partial_{\tau}^{\beta}f_{k}(\tilde x,\tau)\mid\;\leq\;C_{\alpha \beta K} \quad\rm{on }\;R^{N}\times[0,\infty)\;\rm{ for \; any }\;\alpha\;,\;\beta\;,\;K
\end{equation}

By comparing $(\ref{eqn17b})\;, \;(\ref{eqn18b})$ with $(\ref{eqn5a})$ we can see that infinitely differentiable functions $u_{k}^{0}(\tilde x),  f_{k}(\tilde x,\tau)$ 

are satisfying inequalities of the type $(\ref{eqn5a})$ and hence, $u_{k}^{0}(\tilde x) \in S(R^{N}), \; f_{k}(\tilde x,\tau) \in S(R^{N})$.
\nonumber\\

Let us consider the first step of iterative process.

We can see that inner integrals (Fourier transforms) in the integral operators $S_{11}(), S_{12}(), S_{13}(), S_{21}(),$ 

$ S_{22}(), S_{23}(), S_{31}(), S_{32}(), S_{33}(), B()$ from formulas $(\ref{eqn160})\;, (\ref{eqn161})\;, (\ref{eqn162})$ transform $u_{k}^{0}(\tilde x),  f_{k}(\tilde x,\tau)$ into 

$\widehat{u}_{k}^{0}(\gamma) \in S(R^{N}), \;$ $\widehat{f}_{k}(\gamma,\tau) \in S(R^{N})$, 
according to $\cite{gS01}$. 

Multiplication of $\widehat{u}_{k}^{0}(\gamma)$ by $\sz{e} ^{-\nu (\gamma_{1}^{2}+\gamma_{2}^{2} +\gamma_{3}^{2})t}$ and $\widehat{f}_{k}(\gamma,\tau)$ by 
$\sz{e} ^{-\nu (\gamma_{1}^{2} +\gamma_{2}^{2} +\gamma_{3}^{2})(t-\tau)}$ and by fractions 
\nonumber\\

\begin{center}
$\mid\frac{( \gamma_{2}^{2} +\gamma_{3}^{2})} { (\gamma_{1}^{2} +\gamma_{2}^{2}+\gamma_{3}^{2} ) }\mid< 1$,
$\mid\frac{( \gamma_{1}\cdot\gamma_{2})} { (\gamma_{1}^{2} +\gamma_{2}^{2}+\gamma_{3}^{2} ) }\mid< 1$,
$\mid\frac{( \gamma_{1}\cdot\gamma_{3})} { (\gamma_{1}^{2} +\gamma_{2}^{2}+\gamma_{3}^{2} ) }\mid< 1$,
\end{center}
$\nonumber\\$
\begin{center}
$\mid\frac{( \gamma_{3}^{2} +\gamma_{1}^{2})} { (\gamma_{1}^{2} +\gamma_{2}^{2}+\gamma_{3}^{2} ) }\mid< 1$,
$\mid\frac{( \gamma_{2}\cdot\gamma_{3})} { (\gamma_{1}^{2} +\gamma_{2}^{2}+\gamma_{3}^{2} ) }\mid< 1$,
$\mid\frac{( \gamma_{1}^{2} +\gamma_{2}^{2})} { (\gamma_{1}^{2} +\gamma_{2}^{2}+\gamma_{3}^{2} ) }\mid< 1$
\end{center}
$\nonumber\\$

keeps result functions in class $S(R^{N})$. 

Inverse Fourier transforms (outer integrals in the integral-operators $S_{11}(), S_{12}(), S_{13}(), S_{21}(),$

$S_{22}(), S_{23}(), S_{31}(), S_{32}(), S_{33}(), B()$) result in $u_{1k}(x,\tau) \in S(R^{N})$ according 
to $\cite{gS01}$. Integrating $u_{1k}(x,\tau)$ 

with respect to $\tau$ over the interval [0, t] keeps functions $u_{1k}(x,t)$ in class $S(R^{N})$.
\nonumber\\

Let us consider the second step of iterative process.

We obtain the first correction to applied force 

\begin{equation}\label{eqn169a}
\vec{f}_{2}^{*}\;=\;(\vec{u}_{1}\;\cdot\;\nabla)\;\vec{u}_{1}
\end{equation}

from formula $(\ref{eqn169})$. For arbitrary k $({1\leq k \leq N})$ we have 

\begin{equation}\label{eqn7a}
f_{2k}^{*}\;=\;\; \sum_{n=1}^{N} u_{1n}\frac{\partial u_{1k}}{\partial x_{n}}\;\;\;\;\;\;\;\;\;\;\;\; 
\end{equation}

From $u_{1k}(x,\tau) \in S(R^{N})$ it follows that $\frac{\partial u_{1k}}{\partial x_{n}} \in S(R^{N})$, and hence $f_{2k}^{*} \in S(R^{N})$ $\cite{gS01}$.

We can obtain the first correction to velocity 

\begin{equation}\label{eqn170a}
\vec{u}_{2}^{*}\;=\;\bar{\bar{S}}\;\cdot\;\vec{f}_{2}^{*}
\end{equation}

from formula $(\ref{eqn170})$. After analogous reasoning for components $u_{2k}^{*}$ like for $u_{1k}$ on the first step of 

iterative process, we have $u_{2k}^{*}\in S(R^{N})$ according to $\cite{gS01}$.

Hence, we have received that on any arbitrary step $l \;(l > 1)$ of the iterative process a correction to the 

force $\vec{f}_{l}^{*}$ as well as correction to the velocity $\vec{u}_{l}^{*}$ are infinitely differentiable functions
and $f_{lk}^{*}\in S(R^{N}),$ 

$u_{lk}^{*}\in S(R^{N})$ $({1\leq k \leq N})$.
\nonumber\\

\textbf {b)}$\;$ Following $\cite{gS01}$ we introduce classes of functions $W_{M}$ and $W^{\Omega}$  in this paragraph. Let M(x) and $\Omega (t)$
 
be dual functions, in Young's sense, and let  $W_{M}$ be the class of all infinitely differentiable functions

 $\varphi(x)\;(-\infty < x < \infty)$, satisfying the inequalities 

\begin{equation}\label{eqn5b}
\mid \varphi^{(q)}(x)\mid\;\leq\;C_{q}\sz{e} ^{-M(x)}\; ( q = 0, 1, 2, ...).
\end{equation}

where $C_{q}$ is a constant and depends on $\varphi(x)$.

If $\psi(s) = \textbf {F}[\varphi(x)]$ is Fourier transform, then

\begin{equation}\label{eqn5c}
\mid s^{q}\psi(\sigma + i\tau)\mid\;\leq\;C_{q}^{'}\sz{e} ^{\Omega(\tau)}\; ( q = 0, 1, 2, ...).
\end{equation}

Let $W^{\Omega}$ be the class of all entire functions $\psi(s)$ satisfying inequalities of the form $(\ref{eqn5c})$.

Then, $\textbf {F}W_{M} = W^{\Omega}$, in other words, the Fourier transform operator \textbf {F} maps the class $W_{M}$ onto the class 

$W^{\Omega}$ and $\textbf {F}W^{\Omega} = W_{M}$, i.e., the Fourier transform operator F maps the class $W^{\Omega}$ onto the class $W_{M}$ $\cite{gS01}$.
\nonumber\\

Now let us consider $u_{k}^{0}(\tilde x) \in W_{M}(R^{N}), \; f_{k}(\tilde x,\tau) \in W_{M}(R^{N})$ and go to the first step of iterative process. 

Inner integrals (Fourier transforms) in the integral-operators $S_{11}(), S_{12}(), S_{13}(), S_{21}(),$ $ S_{22}(), S_{23}(), S_{31}(), $ 

$S_{32}(), S_{33}(), B()$ from formulas $(\ref{eqn160})\;,$ $(\ref{eqn161})\;, (\ref{eqn162})$ transform $u_{k}^{0}(\tilde x),  f_{k}(\tilde x,\tau)$ into $\widehat{u}_{k}^{0}(\gamma) \in W^{\Omega}(R^{N}), \;$ 

$\widehat{f}_{k}(\gamma,\tau) \in W^{\Omega}(R^{N})$, according to $\cite{gS01}$. 

Multiplication of $\widehat{u}_{k}^{0}(\gamma)$ by $\sz{e} ^{-\nu (\gamma_{1}^{2}+\gamma_{2}^{2} +\gamma_{3}^{2})t}$ and $\widehat{f}_{k}(\gamma,\tau)$ by 
$\sz{e} ^{-\nu (\gamma_{1}^{2} +\gamma_{2}^{2} +\gamma_{3}^{2})(t-\tau)}$ and by fractions 
\nonumber\\

\begin{center}
$\mid\frac{( \gamma_{2}^{2} +\gamma_{3}^{2})} { (\gamma_{1}^{2} +\gamma_{2}^{2}+\gamma_{3}^{2} ) }\mid< 1$,
$\mid\frac{( \gamma_{1}\cdot\gamma_{2})} { (\gamma_{1}^{2} +\gamma_{2}^{2}+\gamma_{3}^{2} ) }\mid< 1$,
$\mid\frac{( \gamma_{1}\cdot\gamma_{3})} { (\gamma_{1}^{2} +\gamma_{2}^{2}+\gamma_{3}^{2} ) }\mid< 1$,
\end{center}
$\nonumber\\$
\begin{center}
$\mid\frac{( \gamma_{3}^{2} +\gamma_{1}^{2})} { (\gamma_{1}^{2} +\gamma_{2}^{2}+\gamma_{3}^{2} ) }\mid< 1$,
$\mid\frac{( \gamma_{2}\cdot\gamma_{3})} { (\gamma_{1}^{2} +\gamma_{2}^{2}+\gamma_{3}^{2} ) }\mid< 1$,
$\mid\frac{( \gamma_{1}^{2} +\gamma_{2}^{2})} { (\gamma_{1}^{2} +\gamma_{2}^{2}+\gamma_{3}^{2} ) }\mid< 1$
\end{center}
$\nonumber\\$

keeps result functions in class $W^{\Omega}(R^{N})$. 

Inverse Fourier transforms (outer integrals in the integral-operators $S_{11}(), S_{12}(), S_{13}(), S_{21}(),$

$S_{22}(), S_{23}(), S_{31}(), S_{32}(), S_{33}(), B()$) result in $u_{1k}(x,\tau) \in W_{M}(R^{N})$ according 
to $\cite{gS01}$. Integrating $u_{1k}(x,\tau)$ 

with respect to $\tau$ over the interval [0, t] keeps functions $u_{1k}(x,t)$ in class $W_{M}(R^{N})$.
\nonumber\\

Let us consider the second step of iterative process.

We obtain the first correction to applied force 

\begin{equation}\label{eqn169c}
\vec{f}_{2}^{*}\;=\;(\vec{u}_{1}\;\cdot\;\nabla)\;\vec{u}_{1}
\end{equation}

from formula $(\ref{eqn169})$. For arbitrary k $({1\leq k \leq N})$ we have 

\begin{equation}\label{eqn7c}
f_{2k}^{*}\;=\;\; \sum_{n=1}^{N} u_{1n}\frac{\partial u_{1k}}{\partial x_{n}}\;\;\;\;\;\;\;\;\;\;\;\; 
\end{equation}

From $u_{1k}(x,\tau) \in W_{M}(R^{N})$ it follows that $\frac{\partial u_{1k}}{\partial x_{n}} \in W_{M}(R^{N})$, and hence $f_{2k}^{*} \in W_{M}(R^{N})$ $\cite{gS01}$.

We can obtain the first correction to velocity 

\begin{equation}\label{eqn170c}
\vec{u}_{2}^{*}\;=\;\bar{\bar{S}}\;\cdot\;\vec{f}_{2}^{*}
\end{equation}

from formula $(\ref{eqn170})$. After analogous reasoning for components $u_{2k}^{*}$ like for $u_{1k}$ on the first step of 

iterative process, we have $u_{2k}^{*}\in W_{M}(R^{N})$ according to $\cite{gS01}$.

Hence, we have received that on any arbitrary step $l \;(l > 1)$ of the iterative process a correction to the 

force $\vec{f}_{l}^{*}$ as well as correction to the velocity $\vec{u}_{l}^{*}$ are infinitely differentiable functions
and $f_{lk}^{*}\in W_{M}(R^{N}),$ 

$u_{lk}^{*}\in W_{M}(R^{N})$ $({1\leq k \leq N})$.
\nonumber\\

\textbf {c)}$\;$Let us estimate superiorly solution of the Cauchy problem for the 3D Navier - Stokes equations

by iterative method. The purposes of this estimation are:

$\;\;\;\;\;$1)	to show convergence of the iterative method;

$\;\;\;\;\;$2)	to obtain analytical form of the first and second steps of the iterative process;

$\;\;\;\;\;$3)	to receive estimated formula for the border of convergence region  of the iterative process in the 

space of system parameters.

We  substitute fractions

\begin{eqnarray}\label{eqn170cc}
\;\;\;\;\;\;\;\;\;\;\mid\frac{( \gamma_{2}^{2} +\gamma_{3}^{2})} { (\gamma_{1}^{2} +\gamma_{2}^{2}+\gamma_{3}^{2} ) }\mid <1,\;\;\;\;\;\;\;
\mid\frac{( \gamma_{1}\cdot\gamma_{2})} { (\gamma_{1}^{2} +\gamma_{2}^{2}+\gamma_{3}^{2} ) }\mid <1,\;\;\;\;\;\;\;
\mid\frac{( \gamma_{1}\cdot\gamma_{3})} { (\gamma_{1}^{2} +\gamma_{2}^{2}+\gamma_{3}^{2} ) }\mid <1,
\nonumber\\
\nonumber\\
\mid\frac{( \gamma_{3}^{2} +\gamma_{1}^{2})} { (\gamma_{1}^{2} +\gamma_{2}^{2}+\gamma_{3}^{2} ) }\mid <1,\;\;\;\;\;\;\;
\mid\frac{( \gamma_{2}\cdot\gamma_{3})} { (\gamma_{1}^{2} +\gamma_{2}^{2}+\gamma_{3}^{2} ) }\mid <1,\;\;\;\;\;\;\;
\mid\frac{( \gamma_{1}^{2} +\gamma_{2}^{2})} { (\gamma_{1}^{2} +\gamma_{2}^{2}+\gamma_{3}^{2} ) }\mid <1
\end{eqnarray}

 by 1 in the integral-operators $S_{11}(), S_{12}(), S_{13}(), S_{21}(),$ $S_{22}(), S_{23}(), S_{31}(), S_{32}(), S_{33}()$ from formulas 

$\;(\ref{eqn160}), \;(\ref{eqn161}),\;(\ref{eqn162})$ for all steps of iterative process. 

Then we take 

\begin{eqnarray}\label{eqn170d}
f_{11}(\tilde x_{1},\tilde x_{2}, \tilde x_{3},\tau)\;=\;F\cdot f(\tau)\cdot\sz{e} ^{-\mu^{2}(\tilde x_{1}^{2} + \tilde x_{2}^{2} + \tilde x_{3}^{2})},\;\;\;\;\;\;F,\mu  - constants, F>0,\mu>0
\nonumber\\
\nonumber\\
f_{12}\;=\;f_{13}\;=\;0,\;\vec{u}_{0}\;=\;0
\quad\quad\quad\quad\quad\quad\quad\quad\quad\quad\quad\quad\quad\quad\quad\quad\quad\quad\quad\quad\quad\quad
\end{eqnarray}

After  Fourier transforms (inner integrals in the integral-operators $S_{11}(), S_{21}(), S_{31}()$ from formulas 

$\;(\ref{eqn160}), \;(\ref{eqn161}),\;(\ref{eqn162})$) we have:

\begin{equation}\label{eqn170f}
\widehat{f}_{11}(\gamma_{1},\gamma_{2}, \gamma_{3},\tau)\;=\;F\cdot f(\tau)\cdot \bigg(\frac{\pi}{\mu^{2}}\bigg)^{3/2} \cdot\sz{e} ^{-\frac{(\gamma_{1}^{2} +  \gamma_{2}^{2} + \gamma_{3}^{2})}{4  \mu^{2}}}
\end{equation}

Now we multiply $\widehat{f}_{11}(\gamma_{1},\gamma_{2}, \gamma_{3},\tau)$ by $\sz{e} ^{-\nu (\gamma_{1}^{2} +\gamma_{2}^{2} +\gamma_{3}^{2})(t-\tau)}$, change order of integration by $\tau$ and $\gamma_1, \gamma_2, \gamma_3$ and after Inverse Fourier transforms (outer integrals in the integral-operators $S_{11}(), S_{21}(), S_{31}()$) we have:

\begin{eqnarray}\label{eqn170g}
\widehat{u}_{11}(x_{1},x_{2},x_{3},\tau)\;=\;\frac{F\cdot f(\tau)}{[4\mu^{2} \nu (t - \tau) + 1]^{3/2}} \cdot\sz{e} ^{\frac{-\mu^{2}(x_{1}^{2} +  x_{2}^{2} + x_{3}^{2})}{[4\mu^{2} \nu (t - \tau) + 1]}},\;\;\; 
\nonumber\\
\nonumber\\
\widehat{u}_{11} = -\widehat{u}_{12},\;\;\; \widehat{u}_{12} = \widehat{u}_{13}
\quad\quad\quad\quad\quad\quad\quad\quad\quad\quad\quad\quad\quad\quad\quad\quad
\end{eqnarray}

Then we get:

\begin{equation}\label{eqn170h}
\widehat{u}_{11}(x_{1},x_{2},x_{3},t)\;=\;F\int_{0}^{t}\frac{f(\tau)}{[4\mu^{2} \nu (t - \tau) + 1]^{3/2}} \cdot\sz{e} ^{\frac{-\mu^{2}(x_{1}^{2} +  x_{2}^{2} + x_{3}^{2})}{[4\mu^{2} \nu (t - \tau) + 1]}}d\tau
\end{equation}

We substitute y for $\tau$: y = $\frac{1}{[4\mu^{2} \nu (t - \tau) + 1]}$, dy = $\frac{4\mu^{2} \nu}{[4\mu^{2} \nu (t - \tau) + 1]^{2}}d\tau$, put $f(y) = y^{1/2}$ and receive after integration:

\begin{eqnarray}\label{eqn170i}
\widehat{u}_{11}(x_{1},x_{2},x_{3},t)\;=\;\frac{F}{4\mu^{4} \nu (x_{1}^{2} +  x_{2}^{2} + x_{3}^{2})} \;\bigg[-\sz{e} ^{-\mu^{2}(x_{1}^{2} +  x_{2}^{2} + x_{3}^{2})} + \sz{e} ^{\frac{-\mu^{2}(x_{1}^{2} +  x_{2}^{2} + x_{3}^{2})}{(4\mu^{2} \nu t +1)}}\bigg]\;=
\nonumber\\
\nonumber\\
=\;\frac{F}{4\mu^{4} \nu (x_{1}^{2} +  x_{2}^{2} + x_{3}^{2})} \;\biggl\{ \gamma[1, \mu^{2}(x_{1}^{2} +  x_{2}^{2} + x_{3}^{2})] - \gamma \big [1, \frac{\mu^{2}(x_{1}^{2} +  x_{2}^{2} + x_{3}^{2})}{(4\mu^{2} \nu t +1)} \big] \biggr\}\;=
\nonumber\\
\nonumber\\
=\;\frac{F}{4\mu^{2} \nu } \;\biggl\{ \Phi[1, 2; -\mu^{2}(x_{1}^{2} +  x_{2}^{2} + x_{3}^{2})] - \frac{1}{(4\mu^{2} \nu t +1)}\cdot\Phi \big [1, 2; \frac{-\mu^{2}(x_{1}^{2} +  x_{2}^{2} + x_{3}^{2})}{(4\mu^{2} \nu t +1)} \big] \biggr\}
\nonumber\\
\nonumber\\
\widehat{u}_{11} = -\widehat{u}_{12},\;\;\; \widehat{u}_{12} = \widehat{u}_{13}
\quad\quad\quad\quad\quad\quad\quad\quad\quad\quad\quad\quad\quad\quad\quad\quad
\end{eqnarray}

Here $\gamma(\alpha, x)$ is the incomplete gamma function $\cite{BE253}$ and $\Phi(a, c; x)$ is a confluent hypergeometric function $\cite{BE153}$.

\begin{eqnarray}\label{eqn171a}
\mid\frac{\partial \widehat{u}_{11}(x_{1},x_{2},x_{3},t)}{\partial x_{n}}\mid\;=\;\frac{F\cdot x_{n}}{2\mu^{4} \nu (x_{1}^{2} +  x_{2}^{2} + x_{3}^{2})}\;\biggl\{ \frac{1}{(x_{1}^{2} +  x_{2}^{2} + x_{3}^{2})}\bigg[-\sz{e} ^{-\mu^{2}(x_{1}^{2} +  x_{2}^{2} + x_{3}^{2})} + \sz{e} ^{\frac{-\mu^{2}(x_{1}^{2} +  x_{2}^{2} + x_{3}^{2})}{(4\mu^{2} \nu t +1)}}\bigg]\;+
\nonumber\\
\nonumber\\
+ \;\mu^{2}\;\bigg[-\sz{e} ^{-\mu^{2}(x_{1}^{2} +  x_{2}^{2} + x_{3}^{2})} + \frac{1}{(4\mu^{2} \nu t +1)}\;\sz{e} ^{\frac{-\mu^{2}(x_{1}^{2} +  x_{2}^{2} + x_{3}^{2})}{(4\mu^{2} \nu t +1)}}\bigg] \biggr \}\;=
\quad\quad\quad\quad\quad\quad\quad\quad\quad\quad
\nonumber\\
\nonumber\\
=\;\frac{F\cdot x_{n}}{4 \nu } \;\biggl\{ \Phi[2, 3; -\mu^{2}(x_{1}^{2} +  x_{2}^{2} + x_{3}^{2})] - \frac{1}{(4\mu^{2} \nu t +1)^{2}}\cdot\Phi \big [2, 3; \frac{-\mu^{2}(x_{1}^{2} +  x_{2}^{2} + x_{3}^{2})}{(4\mu^{2} \nu t +1)} \big] \biggr\}\;\;({1\leq n \leq N})
\end{eqnarray}

Here $\widehat{\vec{u}}_{1}(x_{1},x_{2},x_{3},t)$ and $\frac{\partial\widehat{\vec{u}}_{1}(x_{1},x_{2},x_{3},t)}{\partial x_{n}}$ are estimations from above of velocity and partial derivatives of velocity on the first step of the iterative process. 

Now we do a next level of superior estimation of velocity and partial derivatives of velocity to continue the iterative process. For this goal we receive from formulas $\;(\ref{eqn170i}), \;(\ref{eqn171a})$ with condition $(x_{1}^{2} +  x_{2}^{2} + x_{3}^{2}) > 1$

\begin{equation}\label{eqn171b}
\widetilde{u}_{11}(x_{1},x_{2},x_{3},t)\;=\;\frac{F}{4\mu^{2} \nu } \;\sz{e} ^{\frac{-\mu^{2}(x_{1}^{2} +  x_{2}^{2} + x_{3}^{2})}{(4\mu^{2} \nu t +1)}}
\end{equation}

\begin{equation}\label{eqn171c}
\mid\frac{\partial \widetilde{u}_{11}(x_{1},x_{2},x_{3},t)}{\partial x_{n}}\mid\;=\;\frac{F}{2 \mu^{2}\nu } \; \frac{2}{(4\mu^{2} \nu t +1)} \; \sz{e} ^{\frac{-\mu^{2}(x_{1}^{2} +  x_{2}^{2} + x_{3}^{2})}{(4\mu^{2} \nu t +1)}}
\end{equation}

Then in this case we have from formula $(\ref{eqn7c})$ following estimation from above:

\begin{equation}\label{eqn171d}
\mid\widetilde{f}_{21}^{*}\mid\;=\;\mid\widetilde{f}_{22}^{*}\mid\;=\;\mid\widetilde{f}_{23}^{*}\mid\;=\;\mid \widetilde{u}_{11}\frac{\partial \widetilde{u}_{11}}{\partial x_{n}}\mid\;=\;\frac{F^{2}}{4 \mu^{4}\nu^{2} } \; \frac{1}{(4\mu^{2} \nu t +1)} \; \sz{e} ^{\frac{-2 \mu^{2}(x_{1}^{2} +  x_{2}^{2} + x_{3}^{2})}{(4\mu^{2} \nu t +1)}}
\end{equation}

After  Fourier transforms (inner integrals in the integral-operators $S_{11}(), S_{21}(), S_{31}()$ from formulas 

$\;(\ref{eqn160}), \;(\ref{eqn161}),\;(\ref{eqn162})$) we have:

\begin{eqnarray}\label{eqn170j}
\widetilde{f}_{2k}^{*}(\gamma_{1},\gamma_{2}, \gamma_{3},\tau)\;=\;\frac{F^{2}}{4 \mu^{4}\nu^{2} } \; \frac{1}{(4\mu^{2} \nu \tau +1)} \;\cdot \bigg(\frac{\pi(4\mu^{2} \nu \tau +1)}{2\mu^{2}}\bigg)^{3/2} \cdot\sz{e} ^{-\frac{(4\mu^{2} \nu \tau +1)(\gamma_{1}^{2} +  \gamma_{2}^{2} + \gamma_{3}^{2})}{8  \mu^{2}}}
\nonumber\\
\nonumber\\
({1\leq k \leq N})
\quad\quad\quad\quad\quad\quad\quad\quad\quad\quad\quad\quad\quad\quad\quad\quad
\end{eqnarray}

Now we multiply $\widetilde{f}_{2k}^{*}(\gamma_{1},\gamma_{2}, \gamma_{3},\tau)$ by $\;\sz{e}^{-\nu (\gamma_{1}^{2} +\gamma_{2}^{2} +\gamma_{3}^{2})(t-\tau)}$, apply substitution of fractions $(\ref{eqn170cc})$ by 1, change order of integration by $\tau$ and $\gamma_1, \gamma_2, \gamma_3$ and after Inverse Fourier transforms (outer integrals in the integral-operators $S_{11}(), S_{21}(), S_{31}()$) we have:

\begin{eqnarray}\label{eqn170k}
\widetilde{u}_{2k}^{*}(x_{1},x_{2},x_{3},\tau)\;=\;\frac{F^{2}(4\mu^{2} \nu \tau +1)^{1/2}}{4 \mu^{4}\nu^{2}[8\mu^{2} \nu (t - \tau) + (4\mu^{2} \nu \tau +1)]^{3/2}} \cdot\sz{e} ^{\frac{-2\mu^{2}(x_{1}^{2} +  x_{2}^{2} + x_{3}^{2})}{[8\mu^{2} \nu (t - \tau) + (4\mu^{2} \nu \tau +1)]}},\;\;\; 
\nonumber\\
\nonumber\\
({1\leq k \leq N})
\quad\quad\quad\quad\quad\quad\quad\quad\quad\quad\quad\quad\quad\quad\quad\quad
\end{eqnarray}

Let us further increase level of estimation. To do so we substitute $\sz{e} ^{\frac{-2\mu^{2}(x_{1}^{2} +  x_{2}^{2} + x_{3}^{2})}{[8\mu^{2} \nu (t - \tau) + (4\mu^{2} \nu \tau +1)]}}$ from $(\ref{eqn170k})$ by 
\nonumber\\

$\sz{e} ^{\frac{-2\mu^{2}(x_{1}^{2} +  x_{2}^{2} + x_{3}^{2})}{[8\mu^{2} \nu (t - \tau) + 2(4\mu^{2} \nu \tau +1)]}}\;$= $\;\sz{e} ^{\frac{-\mu^{2}(x_{1}^{2} +  x_{2}^{2} + x_{3}^{2})}{(4\mu^{2} \nu t +1)}}$.

Then we get:

\begin{equation}\label{eqn170l}
\widetilde{u}_{2k}^{*}(x_{1},x_{2},x_{3},t)\;=\;\frac{F^{2}}{4 \mu^{4}\nu^{2}}\cdot\sz{e} ^{\frac{-\mu^{2}(x_{1}^{2} +  x_{2}^{2} + x_{3}^{2})}{(4\mu^{2} \nu t +1)}}\int_{0}^{t}\frac{(4\mu^{2} \nu \tau +1)^{1/2}}{(8\mu^{2} \nu t - 4\mu^{2} \nu \tau +1)^{3/2}} d\tau
\end{equation}

We substitute y for $\tau$: y = $\frac{1}{(8\mu^{2} \nu t - 4\mu^{2} \nu \tau +1)}$, dy = $\frac{4\mu^{2} \nu}{(8\mu^{2} \nu t - 4\mu^{2} \nu \tau +1)^{2}}d\tau$ and receive after integration:

\begin{equation}\label{eqn170z}
\widetilde{u}_{2k}^{*}(x_{1},x_{2},x_{3},t)\;=\;\frac{F^{2}}{8 \mu^{6}\nu^{3}}\cdot\sz{e} ^{\frac{-\mu^{2}(x_{1}^{2} +  x_{2}^{2} + x_{3}^{2})}{(4\mu^{2} \nu t +1)}}\cdot \bigg [1 - \pi/4  
- \frac{1}{(8\mu\nu t + 1)^{1/2}} + arctg \frac{1}{(8\mu\nu t + 1)^{1/2}} \bigg ]
\end{equation}

Now we compare $\widetilde{u}_{11}(x_{1},x_{2},x_{3},t)$ from formula $(\ref{eqn171b})\;$ with $\widetilde{u}_{2k}^{*}(x_{1},x_{2},x_{3},t)$ from formula $(\ref{eqn170z})\;$ and see that the iterative process is converging with estimated condition:

\begin{equation}\label{eqn170x}
\frac{F}{\mu^{4}\nu} < 1
\end{equation}

where F and $\mu$ were introduced in formula $(\ref{eqn170d})$. Condition $(\ref{eqn170x})$ is the estimated formula for the border of convergence region  of the iterative process in the space of system parameters.

For arbitrary step $j$ $(j \geq 2)$ of iterative process we may take $\widetilde{\vec{u}}_{j}$ from formula $(\ref{eqn182})$ and apply estimation algorithm analogous to formulas $(\ref{eqn170cc})$ - $(\ref{eqn170z})$. 
\nonumber\\

Since these results are shown for the superior estimation, then for precise calculations the convergence of the iterative method will be even better.
\nonumber\\

Then we have from formula $(\ref{eqn163})\;$:

\begin{equation}\label{eqn185c}
p_{\infty}\,\;=\; \;\tilde S_{1}(f_{\infty 1})\;+\; \tilde S_{2}(f_{\infty 2})\;+\; \tilde S_{3}(f_{\infty 3})
\end{equation}

Here $\vec{f}_{\infty}$ = ($f_{\infty 1} , f_{\infty 2}, f_{\infty 3}$) is received from formula $(\ref{eqn181a})\;$.

On the other hand we can transform the original system of differential equations $(\ref{eqn7})\; - \;(\ref{eqn9})$ to the equivalent system of integral equations by the scheme of iterative process $(\ref{eqn164})\;, \;(\ref{eqn165})$ for vector $\vec{u}$:

\begin{equation}\label{eqn186}
\vec{u}\;=\;\vec{u}_{1}\;-\;\bar{\bar{S}}\cdot(\vec{u}\;\cdot\;\nabla)\;\vec{u},
\end{equation}

where $\vec{u}_{1}$ is from formula $(\ref{eqn168})$.
We compare the equations $(\ref{eqn185})$ and $(\ref{eqn186})$ and see that the iterative process $(\ref{eqn185})$ converges to the solution of the system $(\ref{eqn186})$ and hence to the solution of the differential equations $(\ref{eqn7})\; - \;(\ref{eqn9})$. 

\textbf{In other words there exist smooth functions}  $\mathbf{p_{\infty}(x, t)}$, $\mathbf{u_{\infty i}(x, t)}$  \textbf{(i = 1, 2, 3) on} $\mathbf{R^{3} \times [0,\infty)}$ \textbf{that satisfy} $\mathbf{(\ref{eqn1}), (\ref{eqn2}), (\ref{eqn3})}$ \textbf{and}

\begin{equation}\label{eqn186b}
\mathbf{p_{\infty}, \;u_{\infty i} \in  C^{\infty}(R^{3} \times [0,\infty)),}
\nonumber\\
\nonumber\\
\end{equation}

\begin{equation}\label{eqn186c}
\mathbf{\int_{R^{3}}|\vec{u}_{\infty}(x, t)|^{2}dx < C } 
\end{equation}

\textbf{for all t} $\mathbf{\geq 0}$.
\nonumber\\

\section{Justification of the analytical iterative method solution for Cauchy problem for the 3D Navier-Stokes equations}\ 
\nonumber\\

\textbf {4.1. Spaces S,}$\;\overrightarrow{\textbf{TS.}}\;\;\cite{GC68}, \cite{RR64}.$ 

Let us consider space S of all infinitely differentiable functions $\varphi$(x) defined in N-dimensional space $R^{N}$ (N = 3), such that when $\;\;\mid x \mid \;\rightarrow\; \infty\;\;$ these functions tend to 0, as well as their derivatives of any order, more rapidly than any power of $\frac{1}{\mid x \mid}$.

To define topology in the space S let us introduce countable system of norms 

\begin{equation}\label{eqn200}
\|\varphi\|_{p}\;=\; \mathop{sup_{\;x}}_{\mid k \mid ,\; \mid q \mid \;\leq \;p}\mid x^{k}D^{q} \varphi(x)\mid\;\;\;(p = 0, 1, 2,...)
\end{equation}

where 

\begin{center}
$\mid x^{k}D^{q} \varphi(x)\mid\; = \; \mid x_{1}^{k_{1}}\ldots x_{N}^{k_{N}}\frac{\partial^{q_{1}+\cdots + q_{N}}\varphi(x)}{\partial {x_{1}^{q_{1}}}\ldots \partial {x_{N}^{q_{N}}} }\mid$
\end{center}
\begin{center}
\end{center}
\begin{center}
$k = (k_{1}, \ldots , k_{N}),\;\; q = (q_{1}, \ldots , q_{N}),\;\;x^{k} =  x_{1}^{k_{1}}\ldots x_{N}^{k_{N}} $
\end{center}
\begin{center}
\end{center}
\begin{center}
$D^{q} = \frac{\partial^{q_{1}+\cdots + q_{N}}}{\partial {x_{1}^{q_{1}}}\ldots \partial {x_{N}^{q_{N}}} },\;\; (k_{1}, \ldots , q_{N} = 0, 1, 2, \ldots) $
\end{center}

Space S is a perfect space (complete countably normed space, in which the bounded sets are compact). Space $\overrightarrow{TS}$ of vector-functions $\vec{\varphi}$ is a direct sum of N  perfect spaces S (N = 3) $\cite{VT80}, \cite{RN72}$ :
\nonumber\\

\begin{center}
$\overrightarrow{TS} = S \oplus S \oplus S$.
\end{center}

To define topology in the space $\overrightarrow{TS}$ let us introduce countable system of norms 

\begin{equation}\label{eqn201}
\|\vec{\varphi}\|_{p}\;=\; \sum_{i = 1}^{N}\|\varphi_{i}\|_{p}\; = \sum_{i = 1}^{N}\mathop{sup_{\;x}}_{\mid k \mid ,\; \mid q \mid \;\leq \;p}\mid x^{k}D^{q} \varphi_{i}(x)\mid\;\;\;(p = 0, 1, 2,...),\;\;(N = 3)
\end{equation}

\textbf {4.2. Equivalence of Cauchy problem in differential form $\textbf {(\ref{eqn1})}$ - $\textbf {(\ref{eqn3})}$ and in the form of an integral equation.}

Let us denote solution of the problem $(\ref{eqn1})$ - $(\ref{eqn3})$ as \{$\vec{u}(x_{1}, x_{2}, x_{3}, t)$, P($x_{1}, x_{2}, x_{3},$ t)\}, in other words let us consider infinitely differentiable by t $\in$ [0,$\infty$) vector-function $\vec{u}(x_{1}, x_{2}, x_{3}, t) \in \overrightarrow{TS}$ and infinitely differentiable function P($x_{1}, x_{2}, x_{3},$ t) $\in$ S, that turn equations $(\ref{eqn1})$ , $(\ref{eqn2})$ into identities. Vector-function $\vec{u}(x_{1}, x_{2}, x_{3}, t)$ also satisfies the initial condition $(\ref{eqn3})\;\; (\vec{u}^{0} (x_{1}, x_{2}, x_{3})\in \overrightarrow{TS})$:

\begin{equation}\label{eqn202}
\vec{u}(x_{1}, x_{2}, x_{3}, t)|_{t = 0}\;=\;\vec{u}^{0}(x_{1}, x_{2}, x_{3})
\end{equation}

Let us put \{$\vec{u}(x_{1}, x_{2}, x_{3}, t)$, P($x_{1}, x_{2}, x_{3},$ t)\} in to equations $(\ref{eqn1})$ , $(\ref{eqn2})$ and apply Fourier and Laplace transforms to the result identities considering initial condition $(\ref{eqn3})$. After all required operations (as in parts 2 and 3) we receive that vector-function $\vec{u}(x_{1}, x_{2}, x_{3}, t)$ satisfies integral equation:

\begin{equation}\label{eqn203}
\vec{u}\;=\;- \bar{\bar{S}}\cdot(\vec{u}\cdot\nabla)\vec{u} + \vec{u}_{1} = \bar{\bar{S}}^{\nabla}\cdot\vec{u}
\end{equation}

where $\vec{u}_{1} \in \overrightarrow{TS}$ is from $(\ref{eqn168})$:

\begin{equation}\label{eqn204}
\vec{u}_{1}\;=\;\bar{\bar{S}}\cdot\vec{f}\;+\;\bar{\bar{B}}\cdot\vec{u}^{0}\;,\;\;\;\;\;\\
\end{equation}

Here $\vec{f} \in \overrightarrow{TS}$, $\vec{u}^{0} \in \overrightarrow{TS}$. Results of operators $\bar{\bar{S}}\cdot\vec{f},\;\; \bar{\bar{B}}\cdot\vec{u}^{0},\;\; \bar{\bar{S}}\cdot(\vec{u}\cdot\nabla)\vec{u}\;$ are also belong $\overrightarrow{TS}$ since Fourier transform maps perfect space $\overrightarrow{TS}$ onto $\overrightarrow{TS}$.
Function P $\in$ S is defined by formula 

\begin{equation}\label{eqn205}
P\;=\;\tilde{\tilde{S}}\cdot\vec{f}\;-\;\tilde{\tilde{S}}\cdot(\vec{u}\cdot\nabla)\vec{u}\;,\;\;\;\;\;\\
\end{equation}

where vector-function $\vec{u}$ is received from $(\ref{eqn203})$.

Here $\bar{\bar{S}},\;\bar{\bar{S}}^{\nabla},\;\bar{\bar{B}},\;\tilde{\tilde{S}}\;$ are matrix integral operators.

Going from the other side, let us assume that $\vec{u}(x_{1}, x_{2}, x_{3}, t)\in \overrightarrow{TS}$ is continuous in t $\in$ [0,$\infty$) solution of integral equation $(\ref{eqn203})$. Integral-operators $S_{ij}\cdot(\vec{u}\cdot\nabla)\vec{u}$ are continuous in t $\in$ [0,$\infty$). From here we receive that according to $(\ref{eqn203})$, $(\ref{eqn204})$

\begin{center}
$\vec{u}(x_{1}, x_{2}, x_{3}, 0)  = \vec{u}^{0}(x_{1}, x_{2}, x_{3})$ 
\nonumber\\
\nonumber\
\end{center}

and also that $\vec{u}(x_{1}, x_{2}, x_{3}, t)$ is differentiable by t $\in$ [0,$\infty$). As described before, the Fourier transform maps perfect space $\overrightarrow{TS}$ on itself. Hence, $\vec{u}(x_{1}, x_{2}, x_{3}, t)$ and $P(x_{1}, x_{2}, x_{3}, t)\;$ from formula $(\ref{eqn205})$ is the solution of the Cauchy problem $(\ref{eqn1})$ - $(\ref{eqn3})$. From here we see that solving the Cauchy problem $(\ref{eqn1})$ - $(\ref{eqn3})$ is equivalent to finding continuous in t $\in$ [0,$\infty$) solution of integral equation $(\ref{eqn203})$ with condition $(\ref{eqn204})$.
\nonumber\\

\textbf {4.3. The fixed point principle.}$\cite{KA64},\;$ $\cite{VT80},\;$ $\cite{WR73},\;$$\; \cite{KS01},\;$$\; \cite{GD03},\;$$\;\cite{ADL97}\;$

Let us use the fixed point principle to prove existence and uniqueness of the solution of integral equation $(\ref{eqn203})$. 

For this purpose we will operate with the following properties of matrix integral operator $\bar{\bar{S}}^{\nabla}$:

1.	Matrix integral operator $\bar{\bar{S}}^{\nabla}$ continuously depends on its parameter t $\in$ [0,$\infty$) (based on formulas $(\ref{eqn160})$ - $(\ref{eqn162})$).

2.	Matrix integral operator $\bar{\bar{S}}^{\nabla}$ maps vector-functions $\vec{u}$ from perfect space $\overrightarrow{TS}$ onto perfect space $\overrightarrow{TS}$. This property directly follows from the properties of Fourier transform $\cite{GC68},\;$ and the form of integrands of integral operators $S_{ij},\; B$ (based on formulas $(\ref{eqn160})$ - $(\ref{eqn162})$ and paragraphs \textbf{a)} and \textbf{b)} in part 3).

3.	Matrix integral operator $\bar{\bar{S}}$ is "quadratic" and therefore we will consider vector-functions $\vec{u}_{1},\; \vec{u}\;$ such that $\mid\vec{u}_{1}\mid \leq 1$, $\;\mid\vec{u}\mid \leq 1$ for any value t $\in$ [0,$\infty$).

For example, the convergence of the iterative process in part 3 is achieved because of condition $(\ref{eqn170x})$:
\nonumber\\
\nonumber\
\begin{center}
$\frac{F}{\mu^{4}\nu} < 1$
\nonumber\\
\nonumber\
\end{center}

$\\\\$
4. $\|\bar{\bar{S}}^{\nabla}\cdot\vec{u}\; - \;\bar{\bar{S}}^{\nabla}\cdot\vec{u}^{'}\|_{p}\;<\;\|\vec{u}\;-\;\vec{u}^{'} \|_{p}\;$ for any $\vec{u},\;\vec{u}^{'} \in \overrightarrow{TS}\;\;$ (${\vec{u}\;\neq \;\vec{u}^{'},\;}$use property 3)  and any t $\in$ [0,$\infty$)    
(based on formulas $(\ref{eqn160})$ - $(\ref{eqn162})$).

Let us prove that matrix integral operator $\bar{\bar{S}}^{\nabla}$ is a contraction operator.
\nonumber\\

\textsc{Theorem 1. }\textbf{Contraction operator.}$\;\;\;\textsc{\cite{KA64}}\;$

Matrix integral operator $\bar{\bar{S}}^{\nabla}$ maps perfect space $\overrightarrow{TS}$ onto perfect space $\overrightarrow{TS}$, and for any $\vec{u},\;\vec{u}^{'}\in \overrightarrow{TS}$ 

(${\vec{u}\;\neq \;\vec{u}^{'}}$)  satisfying property 3, the condition 4 is valid. 

Then matrix integral operator $\bar{\bar{S}}^{\nabla}$ is a contraction operator, i.e. the following condition is true:

\begin{equation}\label{eqn206}
\|\bar{\bar{S}}^{\nabla}\cdot\vec{u}\; -\;\bar{\bar{S}}^{\nabla}\cdot\vec{u}^{'}\|_{p}\;\leq\;\alpha\cdot\|\vec{u}\;-\;\vec{u}^{'} \|_{p}\;
\end{equation}

where $\alpha\;<\;1$ and is independent from $\vec{u},\;\vec{u}^{'} \in \overrightarrow{TS}$ for any t $\in$ [0,$\infty$).
\nonumber\\

Proof by contradiction.

Let us assume that the opposite is true. Then there exist such $\vec{u}_{n},\;\vec{u}^{'}_{n} \in \overrightarrow{TS}$ (n=1,2,$\ldots$), that 

\begin{equation}\label{eqn207}
\|\bar{\bar{S}}^{\nabla}\cdot\vec{u}_{n}\; -\;\bar{\bar{S}}^{\nabla}\cdot\vec{u}^{'}_{n}\|_{p}\;=\;\alpha_{n}\cdot\|\vec{u}_{n}\;-\;\vec{u}^{'}_{n} \|_{p}\;\;\;\;\;\;\;\;(n=1,2,\ldots;\;\;\alpha_{n}\;\rightarrow\;1 )
\end{equation}

Also, because $\overrightarrow{TS}$ is a perfect space, we can consider that $\vec{u}_{n}\;\rightarrow\;\vec{u} \in \overrightarrow{TS}$ and $\vec{u}^{'}_{n}\;\rightarrow\;\vec{u}^{'} \in \overrightarrow{TS}$.
Then 

the limiting result in $(\ref{eqn207})$ would lead to equality 
\nonumber\\
\nonumber\
\begin{center}
$\|\bar{\bar{S}}^{\nabla}\cdot\vec{u}\; - \;\bar{\bar{S}}^{\nabla}\cdot\vec{u}^{'}\|_{p}\;=\;\|\vec{u}\;-\;\vec{u}^{'} \|_{p}\;,$
\nonumber\\
\nonumber\
\end{center}

which is contradicting condition 4. Hence, $\bar{\bar{S}}^{\nabla}$ is a contraction operator. 
\nonumber\\

\textsc{Theorem 2. }\textbf{Existence and uniqueness of solution.}$\;\;\;\textsc{\cite{KA64}}\;$

Let us consider a contraction operator $\bar{\bar{S}}^{\nabla}$. Then there exists a unique solution $\vec{u}^{*}$ of equation $(\ref{eqn203})$ in

space $\overrightarrow{TS}$ for any t $\in$ [0,$\infty$). Also in this case it is possible to obtain $\vec{u}^{*}$ as a limit of sequence  $\{\vec{u}_{n}\}$ , 

where

\begin{center}
$\vec{u}_{n+1}\; =\;\bar{\bar{S}}^{\nabla}(\vec{u}_{n}) \;\;\;\;\;\;\;\;(n=0,1,\ldots),$ 
\nonumber\\
\nonumber\
\end{center}

and $\vec{u}_{0}\; = \; 0$, and $\vec{u}_{1}$ is taken from $(\ref{eqn204})$.

The rate of conversion of the sequence  $\{\vec{u}_{n}\}$  to the solution can be defined from the following inequality:

\begin{equation}\label{eqn208}
\|\vec{u}_{n}\; -\;\vec{u}^{*}\|_{p}\;\leq\;\frac{\alpha^{n}}{(1 - \alpha)}\|\vec{u}_{1}\;-\;\vec{u}_{0} \|_{p}\;\;\;\;\;\;\;\;\;(n=0,1,\ldots)
\end{equation}

Proof:

Since 

\begin{center}
$\vec{u}_{n+1}\; =\;\bar{\bar{S}}^{\nabla}(\vec{u}_{n}) ,\;\;\;\;\vec{u}_{n}\; =\;\bar{\bar{S}}^{\nabla}(\vec{u}_{n-1})$ 
\nonumber\\
\nonumber\
\end{center}

then according to $(\ref{eqn206})$: 
\nonumber\\
\nonumber\
\begin{center}
$\|\vec{u}_{n+1}\; -\;\vec{u}_{n}\|_{p}\;\leq\;\alpha\cdot\|\vec{u}_{n}\;-\;\vec{u}_{n-1} \|_{p}$.
\nonumber\\
\nonumber\
\end{center}

Subsequently using analogous inequalities while decreasing n we will receive: 
\nonumber\\
\nonumber\
\begin{center}
$\|\vec{u}_{n+1}\; -\;\vec{u}_{n}\|_{p}\;\leq\;\alpha^{n}\cdot\|\vec{u}_{1}\;-\;\vec{u}_{0} \|_{p}$.
\nonumber\\
\nonumber\
\end{center}

From this result it follows that

\begin{eqnarray}\label{eqn209}
\|\vec{u}_{n+l}\; -\;\vec{u}_{n}\|_{p}\;\leq\;\|\vec{u}_{n+l}\; -\;\vec{u}_{n+l-1}\|_{p}+\cdots+\|\vec{u}_{n+1}\; -\;\vec{u}_{n}\|_{p}\;\leq\;
\nonumber\\
\nonumber\\
\leq(\alpha^{n+l-1}\;+\cdots+\;\alpha^{n})\|\vec{u}_{1}\;-\;\vec{u}_{0} \|_{p}\leq\frac{\alpha^{n}}{(1 - \alpha)}\|\vec{u}_{1}\;-\;\vec{u}_{0} \|_{p}.
\quad\quad
\end{eqnarray}

Because of $\alpha^{n} \rightarrow 0$ for n$\;\;\rightarrow\; \infty$, the obtained estimation $(\ref{eqn209})$ shows that sequence  $\{\vec{u}_{n}\}$  is a Cauchy 

sequence, and since the space $\overrightarrow{TS}$ is a perfect space, this sequence converges to an element $\vec{u}^{*}\in \overrightarrow{TS}$, 

such that $\bar{\bar{S}}^{\nabla}(\vec{u}^{*})$ makes sense. We use inequality $(\ref{eqn206})$ again and have:
\nonumber\\
\nonumber\
\begin{center}
$\|\vec{u}_{n+1}\; -\;\bar{\bar{S}}^{\nabla}(\vec{u}^{*})\|_{p}\;=\;\|\bar{\bar{S}}^{\nabla}(\vec{u}_{n})\; -\;\bar{\bar{S}}^{\nabla}(\vec{u}^{*})\|_{p}\;\leq\;\alpha\cdot\|\vec{u}_{n}\;-\;\vec{u}^{*} \|_{p}\;\;\;\;\;\;\;\;\;(n=0,1,2,\ldots)$
\nonumber\\
\nonumber\
\end{center}

The right part of the above inequality tends to 0 for n$\;\;\rightarrow\; \infty$, and it means that 
$\;\;\;\vec{u}_{n+1}\; \rightarrow\;\bar{\bar{S}}^{\nabla}(\vec{u}^{*})$

and $\;\;\;\vec{u}^{*}\; =\;\bar{\bar{S}}^{\nabla}(\vec{u}^{*})$. In other words $\;\;\vec{u}^{*}\;\;$ is the solution of equation $(\ref{eqn203})$.

Uniqueness of the solution also follows from $(\ref{eqn206})$.
In fact, if there would exist another solution $\widetilde{\vec{u}} \in \overrightarrow{TS}$, 

then 
\nonumber\\
\nonumber\
\begin{center}
$\|\widetilde{\vec{u}}\; -\;\vec{u}^{*}\|_{p}\;=\;\|\bar{\bar{S}}^{\nabla}(\widetilde{\vec{u}})\; -\;\bar{\bar{S}}^{\nabla}(\vec{u}^{*})\|_{p}\;\leq\;\alpha\cdot\|\widetilde{\vec{u}}\;-\;\vec{u}^{*} \|_{p}.$
\nonumber\\
\nonumber\
\end{center}

Such situation could happen only if $\|\widetilde{\vec{u}}\; -\;\vec{u}^{*}\|_{p}\;=\;0$, or $\widetilde{\vec{u}}\;=\;\vec{u}^{*}$.
\nonumber\\

We can also receive an estimation $(\ref{eqn208})$ from estimation $(\ref{eqn209})$ as a limiting result 
for $l\;\rightarrow\; \infty$.

Now let us show that continuous dependence of operator $\bar{\bar{S}}^{\nabla}$ on t leads to continuous dependence of
the solution of the problem on t.
We will say that matrix integral operator $\bar{\bar{S}}^{\nabla}$ is continuous on t in
the point $t_{0} \in$ [0,$\infty$), if for any sequence $\{t_{n}\} \in$ [0,$\infty$) such as $t_{n}\; \rightarrow\; t_{0}$ for $n\;\rightarrow\; \infty$, the following is true:

\begin{equation}\label{eqn210}
\bar{\bar{S}}_{t_{n}}^{\nabla}(\vec{u})\;\rightarrow\;\bar{\bar{S}}_{t_{0}}^{\nabla}(\vec{u})
\end{equation}

for any $\vec{u} \in \overrightarrow{TS}$.

From \textsc{Theorem 2} it follows that for any t $\in$ [0,$\infty$) equation $(\ref{eqn203})$ has a unique solution, which is depending on t. Let us denote it as $\vec{u}_{t}^{*}$. We will say that solution of equation $(\ref{eqn203})$ is continuously depending on $t$ for $t\;=\;t_{0}$, if for any sequence $\{t_{n}\} \in$ [0,$\infty$); such as $t_{n}\; \rightarrow\; t_{0}$ for $n\;\rightarrow\; \infty$, the following is true:
\nonumber\\
\nonumber\
\begin{center}
$\vec{u}_{t_{n}}^{*}\;\rightarrow\;\vec{u}_{t_{0}}^{*}$.
\nonumber\\
\nonumber\
\end{center}

\textsc{Theorem 3. }\textbf{Continuous dependence of solution on t.}$\;\;\;\textsc{\cite{KA64}}\;$

Let us consider operator $\bar{\bar{S}}_{t}^{\nabla}$ that satisfies condition $(\ref{eqn206})$ for any t $\in$ [0,$\infty$), with $\alpha$ independent on t 

and that operator $\bar{\bar{S}}_{t}^{\nabla}$ is continuous on t in a point $t_{0} \in$ [0,$\infty$).
Then for $t\;=\;t_{0}$ the solution of equation 

$(\ref{eqn203})$ is continuously depending on t.

Proof: 

Let us consider any t $\in$ [0,$\infty$). We will construct the solution $\vec{u}_{t}^{*}$ of equation $(\ref{eqn203})$ as a limit of sequence 

 $\{\vec{u}_{n}\}$ :

\begin{equation}\label{eqn211}
\vec{u}_{n+1}\; =\;\bar{\bar{S}}_{t}^{\nabla}(\vec{u}_{n}) \;\;\;\;\;\;\;\;(n=0,1,\ldots;\;\;\;\vec{u}_{0}\;=\;\vec{u}_{t_{0}}^{*})
\end{equation}

Let us rewrite inequality $(\ref{eqn208})$ for n = 0:

\begin{equation}\label{eqn212}
\|\vec{u}^{*}\; -\;\vec{u}_{0}\|_{p}\;\leq\;\frac{1}{(1 - \alpha)}\|\vec{u}_{1}\;-\;\vec{u}_{0} \|_{p}
\end{equation}

Since $\vec{u}_{t_{0}}^{*}\; =\;\bar{\bar{S}}_{t_{0}}^{\nabla}(\vec{u}_{t_{0}}^{*})$, then because of $(\ref{eqn211})$ and $(\ref{eqn212})$ we have:

\begin{equation}\label{eqn213}
\|\vec{u}_{t}^{*}\; -\;\vec{u}_{t_{0}}^{*}\|_{p}\;\leq\;\frac{1}{(1 - \alpha)}\|\vec{u}_{1}\;-\;\vec{u}_{0} \|_{p}\;=\;\frac{1}{(1 - \alpha)}\|\bar{\bar{S}}_{t}^{\nabla}(\vec{u}_{t_{0}}^{*})\;-\;\bar{\bar{S}}_{t_{0}}^{\nabla}(\vec{u}_{t_{0}}^{*}) \|_{p}
\end{equation}

Now with the help of $(\ref{eqn210})$ we obtain the required continuity of $\vec{u}_{t}$ for $t\;=\;t_{0}$.
\nonumber\\

Let us consider now a norm $L_{2}$ in space $\overrightarrow{TS}$:
\nonumber\\
\nonumber\
\begin{center}
$\|\vec{u}(x, t)\|_{L_{2}}\;=\;\bigg(\int_{R^{3}}|\vec{u}(x, t)|^{2}dx \bigg)^{1/2}$.
\nonumber\\
\nonumber\
\end{center}

Then we can rewrite estimation $(\ref{eqn186c})$ in the following form:
\nonumber\\
\nonumber\
\begin{center}
$\|\vec{u}(x, t)\|_{L_{2}}^{2}\;=\;\bigg(\int_{R^{3}}|\vec{u}(x, t)|^{2}dx \bigg)\;<\;C$.
\nonumber\\
\nonumber\
\end{center}

From the other side $\|\cdot\|_{L_{2}}$ is a weaker norm than $\|\cdot\|_{p}$. 

Solution of the problem exists in norm $\|\cdot\|_{p}$, which means that $\|\vec{u}(x, t)\|_{p}\;<\;\infty$. Hence there exists solution in norm $\|\cdot\|_{L_{2}}$, and hence condition $(\ref{eqn186c})$ 
\nonumber\\
\nonumber\
\begin{center}
$\bigg(\int_{R^{3}}|\vec{u}(x, t)|^{2}dx \bigg)\;<\;C$
\nonumber\\
\nonumber\
\end{center}

is true for any t $\in$ [0,$\infty$).

\appendix\section*{}

The Fourier integral can be stated in the forms: 

\begin{eqnarray}\label{A3}
U( \gamma_{1} , \gamma_{2} , \gamma_{3})=F[\, u(x_{1} , x_{2} , x_{3})]= \frac{1}{(2\pi)^{3/2}} \int_{-\infty}^{\infty} \int_{-\infty}^{\infty} \int_{-\infty}^{\infty} u( x_{1} , x_{2} , x_{3})\,\sz{e}  ^{ i( \gamma_{1} x_{1} + \gamma_{2} x_{2} + \gamma_{3} x_{3}) } dx_{1} dx_{2} dx_{3} 
\nonumber\\
\nonumber\\
u( x_{1} , x_{2} , x_{3})= \frac{1}{(2\pi)^{3/2}} \int_{-\infty}^{\infty} \int_{-\infty}^{\infty} \int_{-\infty}^{\infty} U( \gamma_{1} , \gamma_{2},\gamma_{3}  )\, \sz{e}  ^{- i( \gamma_{1} x_{1} + \gamma_{2} x_{2} + \gamma_{3} x_{3}) } d\gamma_{1} d\gamma_{2} d\gamma_{3} 
\end{eqnarray}
$\\\\$
The Laplace integral is usually stated in the following form:

\begin{equation}\label{A4}
U^{\otimes}(\eta)=L[\,u(t)\,]= \int_{0}^{\infty}u(t)\, \sz{e}  ^{-\eta t}dt
\;\;\;\;\; u(t)=\frac{1}{2\pi i}\int_{c- i \infty }^{c + i \infty} U^{\otimes}(\eta) \,\sz{e}  ^{\eta t}d\eta \;\;\;\;\; c > c_{0}
\end{equation}

\begin{equation}\label{A5}
L[\,u^{'}(t)\,]=\eta \,U^{\otimes}(\eta)-u(0)
\end{equation}

$\\\\$\textsc{The convolution theorem A.1.}$\;\;\;\textsc{{\cite{DP65}}},\;\textsc{{\cite{DW46}}}$
$\\\\$
If integrals
\[ U_{1}^{\otimes}(\eta)= \int_{0}^{\infty}u_{1}(t)\, \sz{e}  ^{-\eta t}d\,t  \;\;\;\;\;\;\;\;\;\; U_{2}^{\otimes}(\eta)= \int_{0}^{\infty}u_{2}(t)\, \sz{e}  ^{-\eta t}d\,t \]

absolutely converge by $Re\, \eta > \sigma_{d}$, then  $U^{\otimes}(\eta)\,= \,U_{1}^{\otimes}(\eta)\, U_{2}^{\otimes}(\eta)$ is Laplace transform of 

\begin{equation}\label{A6}
u(t)=\int_{0}^{t}u_{1}(t-\tau)\,u_{2}(\tau)\,d\,\tau
\\
\end{equation}

Useful \emph{Laplace integral}:

\begin{equation}\label{A7}
L[\,\sz{e}  ^{\eta_{k}t}\,]\,=\,\int_{0}^{\infty}\sz{e}  ^{-(\eta-\eta_{k})\,t}d\,t
\;=\; \frac{1}{(\eta-\eta_{k})}\;\;\;\;\;\;\;\;\;(Re\,\eta\,>\,\eta_{k})
\end{equation}
\nonumber\\

\textbf{Acknowledgment}:  We express our sincere gratitude to Professor L. Nirenberg, whose suggestion led to conduction of this research.
We are also very thankful to Professor P.G. Lemarié-Rieusset for his valuable comments to our work.

\end{document}